%% file: vertexcuts.tex
\documentclass{article}
\usepackage{amsmath}
\usepackage{amssymb}
\usepackage{amsfonts}
\usepackage{amsthm}
\usepackage[retainorgcmds]{IEEEtrantools}
\usepackage{graphicx}
\usepackage{caption}
\usepackage{subcaption}
\DeclareGraphicsExtensions{.eps}
\theoremstyle{definition}
\newtheorem{preldef1}{Definition}[section]
\newtheorem{preldef2}[preldef1]{Definition}
\newtheorem{preldef3}[preldef1]{Definition}
\newtheorem{crossdef}{Definition}[section]
\newtheorem{halfdefn}{Definition}[section]
\newtheorem{halfdef2}[halfdefn]{Definition}
\newtheorem{sepdef}{Definition}[section]
\newtheorem{sepdef2}[sepdef]{Definition}
\newtheorem{ringdef}{Definition}[section]
\newtheorem{nvertexdef}[ringdef]{Definition}
\newtheorem{pretdef1}{Definition}[section]
\newtheorem{pretdef2}[pretdef1]{Definition}
\newtheorem{succdef}{Definition}[section]
\theoremstyle{plain}
\newtheorem{prel1}[preldef1]{Lemma}
\newtheorem{prel2}[preldef1]{Lemma}
\newtheorem{prel3}[preldef1]{Lemma}
\newtheorem{prel4}[preldef1]{Lemma}
\newtheorem{prel5}[preldef1]{Lemma}
\newtheorem{prel6}[preldef1]{Lemma}
\newtheorem{prel7}[preldef1]{Lemma}
\newtheorem{Bcuts}[crossdef]{Lemma}
\newtheorem{cross1}[crossdef]{Lemma}
\newtheorem{cross2}[crossdef]{Lemma}
\newtheorem{cross3}[crossdef]{Lemma}
\newtheorem{cross4}[crossdef]{Lemma}
\newtheorem{half1}[halfdefn]{Lemma}
\newtheorem{half2}[halfdefn]{Corollary}
\newtheorem{half3}[halfdefn]{Lemma}
\newtheorem{half4}[halfdefn]{Lemma}
\newtheorem{half5}[halfdefn]{Lemma}
\newtheorem{half6}[halfdefn]{Lemma}
\newtheorem{half7}[halfdefn]{Lemma}
\newtheorem{sep1}[sepdef]{Lemma}
\newtheorem{sep2}[sepdef]{Lemma}
\newtheorem{sep3}[sepdef]{Lemma}
\newtheorem{sep4}[sepdef]{Lemma}
\newtheorem{sep5}[sepdef]{Proposition}
\newtheorem{pret1}[pretdef1]{Lemma}
\newtheorem{pret2}[pretdef1]{Lemma}
\newtheorem{pret3}[pretdef1]{Lemma}
\newtheorem{pret4}[pretdef1]{Lemma}
\newtheorem{pret5}[pretdef1]{Lemma}
\newtheorem*{menger}{Menger's Theorem}
\newtheorem{pret6}[pretdef1]{Lemma}
\newtheorem{succthm}[succdef]{Theorem}
\newtheorem{succthm2}[succdef]{Theorem}
\newtheorem{succthm3} [succdef]{Theorem}
\newtheorem*{succulent}{Theorem \ref{sucthm}}
\newtheorem*{stallings}{Stallings's Theorem}
\newtheorem{cactends}{Theorem}[section]
\newtheorem{cactedge}[cactends]{Theorem}
\newcommand{\hashc}{\#-class\ }
\newcommand{\qe}{quasi-equivalent}
\newcommand{\pn}[2]{$#1^{(#2)}$}
\newcommand{\pnm}[2]{#1^{(#2)}}
\newcommand{\xysep}{tight $x$-$y$-separator }
\newcommand{\xyseps}{tight $x$-$y$-separators }
\begin{document}
\title{On the structure of vertex cuts separating the ends of a graph}
\author{Gareth R. Wilkes}
\date{October 10, 2013}
\maketitle
\begin{abstract}
Dinits, Karzanov and Lomonosov showed that the minimal edge cuts of a finite graph have the structure of a cactus, a tree-like graph constructed from cycles. Evangelidou and Papasoglu extended this to minimal cuts separating the ends of an infinite graph. In this paper we show that minimal vertex cuts separating the ends of a graph can be encoded by a succulent, a mild generalization of a cactus that is still tree-like. We go on to show that the earlier cactus results follow from our work.
\end{abstract}
\section{Introduction and Definitions}
\input{intro}
\section{Preliminaries}
\input{prelims2}
\section{Crossing Cuts}
\input{crosscuts}
\section{Half-Cuts}
\input{halfcuts}
\section{Separation Systems}
\input{separations}
\section{The structure of a \hashc}
\input{classstructure}
\section{Pretrees}
\input{pretrees}
\section{Succulents}
\input{succulents}
\section{Applications}
\input{applications}

\bibliographystyle{plain}
\bibliography{vertexcuts}

\begin{center}
\emph{Email address}: \texttt{gareth.wilkes@sjc.ox.ac.uk} \\
\end{center}
\end{document}

%% file: intro.tex
Lying on the boundaries of several topic areas, vertex and edge cuts of graphs have been considered by graph theorists, network theorists, topologists and geometric group theorists and the study of their structure has led to applications ranging from algorithms to classical group theoretic propositions.\\
Vertex cut pairs of finite graphs were studied by Tutte \cite{tutte1984graph}, who showed that graphs possessing such cuts can be modelled with a tree. This was extended to infinite, locally finite graphs in \cite{droms1995structure}. Dunwoody and Kr\"on \cite{dunwoody2009vertex} then extended this work to cuts of other cardinalities, using vertex cuts to associate structure trees to graphs in a more general context.\\
This process of finding trees associated to graphs gives a way into geometric group theory. If for instance we find a structure tree for the Cayley graph of a group, then in light of the work of Bass and Serre \cite{serre1980trees} we can obtain information about the group from its action on the tree. An example is Stallings's theorem \cite{stallings1968torsion} on the classification of groups with many ends. The work of Dunwoody and Kr\"on \cite{dunwoody2009vertex} and of Papasoglu and Evangelidou \cite{papaedgecuts} yields more proofs of Stallings's theorem along these lines.\\
Dinits, Karzanov and Lomonosov \cite{dinits1991structure} showed that minimal edge cuts of a finite graph have, in addition to a tree-like nature, the finer structure of a cactus graph. For a recent elementary proof see \cite{fleiner2009quick}. Papasoglu and Evangelidou \cite{papaedgecuts} extended this, encoding all minimal edge cuts separating the ends of an infinite graph by a cactus. The important stages in these proofs involve showing that certain collections of `crossing' cuts have a circular structure.\\
In this paper we switch our attention to vertex cuts, showing that we can encode all minimal vertex end cuts of a graph by a tree-like structure called a succulent, which is a mild generalisation of a cactus. A traditional cactus is composed of cycles joined together at vertices in a tree-like fashion. For our succulents we also allow cycles to attach along a single edge, again in a tree-like way. Once again the key step is to show that crossing cuts have a cyclic nature.\\
We will also show how the earlier cactus theorems can be regarded as special cases of our work, and discuss an application to certain finite graphs.\\

Let $\Gamma = (V,E)$ be a connected graph. If $K \subseteq V$ is a set of vertices of the graph, we denote by $\Gamma - K$ the graph obtained from $\Gamma$ by removing $K$ and all edges incident to $K$. $K$ is called a \emph{vertex cut} if $K$ is finite and $\Gamma - K$ is not connected. If $A,B\subseteq\Gamma$ then we say $K$ \emph{separates} $A,B$ if any path joining a vertex of $A$ to a vertex of $B$ intersects $K$.\\
A \emph{ray} of $\Gamma$ is an infinite sequence of distinct consecutive vertices of $\Gamma$. We say that two rays $r_1, r_2$ are \emph{equivalent} if for any vertex cut $K$ all vertices of $r_1 \cup r_2$ except finitely many are contained in the same component of $\Gamma - K$. The \emph{ends} of $\Gamma$ are equivalence classes of rays. If $K$ is a vertex cut of $\Gamma$, we say $K$ is an \emph{end cut} if there are at least two components of $\Gamma - K$ which contain rays.  We say that an end cut is a \emph{mincut} if its cardinality is minimal amongst end cuts of $\Gamma$. A mincut is said to \emph{separate} ends $e_1, e_2$ of a graph if there are rays $r_1, r_2$ representing $e_1, e_2$ respectively such that $r_1, r_2$ are contained in different components of $\Gamma- K$. A mincut gives a partition of the set $\mathcal{E}$ of ends of the graph. Two mincuts are called \emph{equivalent} if they give the same partition of $\mathcal{E}$. We denote the equivalence class of a mincut $K$ by $[K]$, and write $K\sim L$ if $K,L$ are equivalent.\\
A \emph{succulent} is a graph constructed from cycles by joining cycles together at vertices or at single edges, in a `tree-like' fashion. We give a more formal definition of this as definition \ref{sucdef} below. An \emph{end vertex} of a succulent is one joined to at most two edges. We now state the main theorem of the paper.
\begin{succulent}
Let $\Gamma$ be a connected graph such that there are vertex end cuts of $\Gamma$ with finite cardinality. There is a succulent $\mathcal{S}$ with the following properties:
\begin{enumerate}
\item There is a subset $A$ of vertices of $\mathcal{S}$ called the \emph{anchors} of $\mathcal{S}$. If two anchors are adjacent, one of them is an end vertex of the graph. Every vertex of $\mathcal{S}$ not in $A$ is adjacent to an anchor. We define an anchor cut of $\mathcal{S}$ to be a vertex cut containing no anchors which separates some anchors of $\mathcal{S}$. We say anchor cuts are equivalent if they partition $A$ in the same way.
\item There is an onto map $f$ from the ends of $\Gamma$ to the union of the ends of $\mathcal{S}$ with the end vertices of $\mathcal{S}$ which are anchors.
\item There is a bijective map $g$ from equivalence classes of minimal end cuts of $\Gamma$ to equivalence classes of minimal anchor cuts of $\mathcal{S}$ such that ends $e_1, e_2$ of $\Gamma$ are separated by $[K]$ if and only if $f(e_1), f(e_2)$ are separated by $g([K])$.
\item Any automorphism of $\Gamma$ induces an automorphism of $\mathcal{S}$.
\end{enumerate}
\end{succulent}
The author would like to thank Panos Papasoglu for suggesting this problem, and Jonathan Bradford and Leo Wright for proof-reading the paper. The author is also grateful to the EPSRC for funding this research.

%% file: prelims2.tex
\begin{preldef1}
Given a mincut $K$, we call a component of $\Gamma - K$ \emph{proper} if it contains an end, and a \emph{slice} if not. Given a set of vertices $C$, its \emph{boundary} $\partial C$ is the set of those vertices not in $C$ but adjacent to a vertex of $C$; and $C^*=V(\Gamma) - (C\cup \partial C)$.
\end{preldef1}
It will be convenient to assume our graph contains no slices. In the following lemmas we show that we can do this by replacing $\Gamma$ with another graph $\hat{\Gamma}$ which has the same ends and cuts, but no slices. The results in this section are adapted for our needs from more general results proved by Dunwoody and Kr\"on \cite{dunwoody2009vertex}.
\begin{prel1}\label{prel1}
Let $K$, $L$ be mincuts and $C,D$ proper components of $\Gamma - K$, $\Gamma - L$. Suppose that both $C\cap D$ and $C^*\cap D^*$ contain an end. Then $\partial (C\cap D)$, $\partial (C^*\cap D^*)$ are mincuts, 
\begin{equation*}
\partial (C\cap D) = (C\cap L)\cup (K\cap L)\cup(K\cap D)
\end{equation*}
\begin{equation*}
\partial (C^*\cap D^*) = (C^*\cap L)\cup (K\cap L)\cup(K\cap D^*)
\end{equation*}
and 
\begin{eqnarray*}
|C\cap L| = |K\cap D^*|\\
|D\cap K| = |L\cap C^*|
\end{eqnarray*}
\begin{proof}
The boundaries $\partial (C\cap D), \partial (C^*\cap D^*)$ are certainly end cuts, with 
\begin{equation*}
\partial (C\cap D) \subseteq (C\cap L)\cup (K\cap L)\cup(K\cap D)
\end{equation*}
\begin{equation*}
\partial (C^*\cap D^*) \subseteq (C^*\cap L)\cup (K\cap L)\cup(K\cap D^*)
\end{equation*}
Consider the following diagram, where $a, b, c, d, u$ denote the cardinalities of the indicated sets. Let $n$ be the cardinality of a mincut.
\begin{figure}[htp]
\centering
\includegraphics[width=0.3\textwidth]{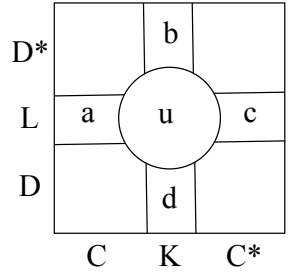}
\caption{Diagram for lemma \ref{prel1}}
\end{figure}
Then
\begin{eqnarray*}
a + c + u = n\\
b + d + u = n
\end{eqnarray*}
Since $(C\cap L)\cup (K\cap L)\cup(K\cap D)$ is an end cut, we have
\begin{equation*}
d + a + u\geq n
\end{equation*}
and similarly
\begin{equation*}
b + c + u\geq n
\end{equation*}
Summing these and comparing with the equalities above, we find them to be equalities; it follows that $a=b, c=d$.
\end{proof}
\end{prel1}
An analogous result holds when $C^*\cap D$, $D^*\cap C$ both contain ends.\\
\begin{prel2}\label{DK36}
If $C, D$ are proper components of cuts $K,L$ then there is a proper component of $\Gamma - K$ containing $C^*\cap L$.
\begin{proof}
Since $K,L$ are end cuts, one of the pairs $\{C\cap D, C^*\cap D^*\}$, $\{C^*\cap D, C\cap D^*\}$ both contain ends. Let $A$ be the appropriate one of $C^*\cap D^*$, $C^*\cap D$. Then using \ref{prel1}, $\partial A$ is a mincut. Let $E$ be a component of $A$ containing an end; $\partial E=\partial A$ is a mincut. Let $C^*_0$ be the component of $C^*$ containing $E$. By \ref{prel1} every vertex $x\in C^*\cap L$ is adjacent to $E$, so $x\in C^*_0$ and $C^*\cap L\subseteq C^*_0$.
\end{proof}
\end{prel2}
\begin{prel3}\label{DK39}
A slice component of a mincut has empty intersection with each mincut. Distinct slices are disjoint. If $Q$ is a slice, then no pair of elements of $\partial Q$ are separated by any mincut.
\begin{proof}
Let $Q_1$ be a slice component of $\Gamma - K$ for a mincut $K$ and let $L$ be a mincut, with a proper component $D$ of $\Gamma - L$. By \ref{DK36} there is a proper component of $\Gamma - K$ containing $C^*\cap L$, and $C$ a proper component containes $C\cap L$. $Q$ is disjoint from both of these, so $Q_1\cap L=\emptyset$.\\
Suppose $Q_2$ is a slice component of $\Gamma - L$. We have $\partial Q_2 \subseteq L$, $\partial Q_1\subseteq K$, hence $Q_1\cap\partial Q_2$, $Q_2\cap\partial Q_1$ are both empty. The components $Q_1, Q_2$ are connected, so this implies that they are disjoint or equal.\\
Finally suppose $x,y\in\partial Q$ for a slice component $Q$ of $\Gamma - K$ and $x,y$ are separated by a mincut $L$. The slice $Q$ is connected so there is a path in $Q$ from $x$ to $y$, which must intersect $L$, but we have seen this is impossible.
\end{proof}
\end{prel3}
We will now show how to replace $\Gamma$ with another graph $\hat{\Gamma}$ which has the same ends and cuts, but no slices. The vertex set $\hat{V}$ of $\hat{\Gamma}$ consists of those vertices of $\Gamma$ which are contained in no slice. Two vertices $u,v\in\hat{V}$ are joined by an edge in $\hat{\Gamma}$ iff they are joined by an edge in $\Gamma$ or if $u,v$ lie in the boundary of some slice of $\Gamma$.
\begin{prel4}
The graph $\hat{\Gamma}$ is connected and the mincuts of $\hat{\Gamma}$ are the same as the mincuts of $\Gamma$. There are no slices in $\hat{\Gamma}$. The ends of $\hat{\Gamma}$ are in bijection with the ends of $\Gamma$.
\begin{proof}
First we show that if $K$ is a mincut and $C$ is a proper component thereof, then $\partial \hat{C}$, the boundary of $\hat{C}=C\cap\hat{\Gamma}$ as a subset of $\hat{\Gamma}$, is equal to $K$.\\
Suppose there is $x\in \partial \hat{C} - K$. If $x\in C$ then $x\in\hat{C}$, so $x\in C^*$. Also $x\in\partial\hat{C}$ so there is $y\in\hat{C}$ adjacent to $x$ in $\hat{\Gamma}$. Then there is an edge from $x$ to $y$ in $\hat{\Gamma}$, but not in $\Gamma$; so $x,y$ lie in the boundary of some slice $Q$ of $\Gamma$. By \ref{DK39}, $K\cap Q=\emptyset$. The slice $Q$ is connected, and $Q$ intersects $C$ (at $y$), so $Q\subseteq C$. We then have a path from $x$ to $y$ in $\Gamma$ which is contained in $Q$ except for its endpoint $x$, which is a path from $C^*$ to $C$ not intersecting $K$, a contradiction.\\
Suppose $\exists x\in K - \partial\hat{C}$; $x$ has a neighbour $y$ in $C - \hat{C}$. Then $y$ is contained in a slice component $Q$ of $\Gamma - L$ for a mincut $L$. If $K=L$ then $C, Q$ are disjoint; but $y\in C\cap Q$. So $K\neq L$ and since $Q\subseteq C$ ($Q$ does not intersect $K$ but does intersect $C$) there is $z\in C\cap\partial Q\subseteq C\cap L$. $z$ is not in any slice, so $z\in\hat{C}$. Then $x,z$ are adjacent in $\hat{\Gamma}$; but $z\in \hat{C}$, $x\notin \hat{C}\cup\partial\hat{C}$. Contradiction.\\
Let us discuss the ends of $\hat{\Gamma}$. By definition, slices contain no rays. Thus if $r$ is any ray in $\Gamma$, we can form a new ray in $\hat{\Gamma}$ by deleting any vertices in a slice; the extra edges added in the construction of $\hat{\Gamma}$ will ensure that this is a bona fide ray. If two rays are separated by a (not necessarily minimal) end cut $K$ in $\Gamma$, then the union of $K\cap V(\hat{\Gamma})$ with the boundaries of any slices intersecting $K$ gives an end cut separating the images of the rays in $\hat{\Gamma}$. Similarly, if two rays in $\hat{\Gamma}$ are separated by an end cut $K$ in $\hat{\Gamma}$, then taking the union of $K$ with any slice boundaries intersecting $K$ gives an end cut separating the same rays in $\Gamma$. It follows that the ends of $\hat{\Gamma}$ are in a natural bijection with those of $\Gamma$.\\
The end cuts of $\hat{\Gamma}$ inherited from mincuts of $\Gamma$ are indeed the minimal end cuts of $\hat{\Gamma}$. Suppose $K$ is an end cut of $\hat{\Gamma}$ which is not also an end cut of $\Gamma$. Then two proper components of $\hat{\Gamma}-K$ are connected in $\Gamma$. A path between them can only not intersect $K$ if it passes through a slice $Q$; but points on the boundary of $Q$ are connected in $\hat{\Gamma}$ so we get a path between the two components in $\hat{\Gamma}$ as well. Contradiction. So all mincuts of $\hat{\Gamma}$ are mincuts of $\Gamma$ as well, so the notion of minimality carries over to $\hat{\Gamma}$ too.\\
Finally, there are no slices in $\hat{\Gamma}$. Let $C$ be a component of $\hat{\Gamma}-K$ for a mincut $K$ of $\hat{\Gamma}$ (equivalently of $\Gamma$). Let $C'$ be the component of $\Gamma - K$ containing $C$. $C'$ cannot be a slice as it intersects $V(\hat{\Gamma})$. So $C'$ contains an end of $\Gamma$, whence from above so does $C$. So $C$ is not a slice. 
\end{proof}
\end{prel4}
For the rest of the paper we replace $\Gamma$ with $\hat{\Gamma}$. As we have seen, the ends and cuts of the two graphs are the same, and this is all the structure with which we are concerned, so we lose nothing by doing this. All components of a cut are now proper.\\
We now start to prove some basic properties of mincuts, putting restrictions on cuts which `interact' with each other in some sense, and showing that a mincut does not `interact' with any but finitely many other mincuts. We first define what it means for cuts to not `interact' with each other. We are still following Dunwoody and Kr\"on \cite{dunwoody2009vertex} here, with some minor modifications.
\begin{preldef2}\label{nestdef}
Two cuts $K,L$ are called \emph{nested} if there are components $E,F$ of $\Gamma - K$, $\Gamma-L$ respectively with $E\subseteq F$ or $F\subseteq E$.
\end{preldef2}
Note that if $K,L$ are nested and not equal with say $E\subseteq F$ then all components of $\Gamma - L$ except $F$ are contained in the same component of $\Gamma - K$. This follows since there is an element of $L$ in $E^*$, and by minimality all components of $\Gamma - L$ except $F$ are connected to this vertex by paths which do not intersect $F\cup L$, hence do not intersect $E\cup K$. Note also that these components are still connected in $\Gamma - K$ by similar reasoning. Conversely, all components of $\Gamma - K$ except one are contained in $F$.
\begin{preldef3}
A mincut is called an \emph{A-cut} if it is nested with all other mincuts. It is called a \emph{B-cut} if it separates $\Gamma$ into exactly two components. 
\end{preldef3}
\begin{prel5}\label{ABthm}
A mincut is either an A-cut or a B-cut.
\begin{proof}
Let $K$ be a mincut which is not an A-cut. Then there is a mincut $L$ with which $K$ is not nested. Let $C$ be a (proper) component of $\Gamma - K$, $D$ a (proper) component of $\Gamma - L$. By \ref{DK36} there is a component $C^*_0$ of $C^*$ containing $C^*\cap L$. We wish to show this is the only component of $C^*$. If there is another one $C^*_1$ then $C^*_1\cap L$ is empty; $C^*_1$ is connected so $C^*_1\subseteq D$ or $C^*_1\subseteq D^*$. In the first case, $K,L$ are nested; so the second one happens whichever component $C^*_1$ we choose. So $D\cap C^*\subseteq C^*_0$. Also, $K=\partial C^*_1$ by minimality, so $C^*_1\subseteq D^*$ implies $K\cap D=\emptyset$. Then $D\subseteq C$ or $D\subseteq C^*$ (whence $D\subseteq C^*_0$), in either case $K$ and $L$ are nested. This is a contradiction, so $K$ is a B-cut.
\end{proof}
\end{prel5}
We call a set $S$ of vertices a \emph{tight $x$-$y$-separator} if $\Gamma - S$ has two distinct components $A,B$ which are adjacent to all elements of $S$, with $x\in A, y\in B$.
\begin{prel6}
For each integer $k$ and every pair $x,y$ of vertices of a graph, there are only finitely many \xyseps of order $k$.
\begin{proof}
We proceed by induction. If we take a path from $x$ to $y$, any \xysep of order 1 would have to be a vertex on this path, so there are only finitely many of these.\\
Suppose the lemma holds for all \xyseps of order $k$ in all connected graphs. Take a path $\pi$ from $x$ to $y$ in a graph $\Gamma$ and suppose there are infinitely many \xyseps of order $k+1\geq 2$. Then there is a vertex $z\in \pi - \{x,y\}$ which is contained in infinitely many of these separators. If $S_1, S_2$ are distinct such \xyseps of order $k+1$ in $\Gamma$ then $S_1 - \{z\}$, $S_2 - \{z\}$ are distinct \xyseps of order $k$ in $\Gamma - \{z\}$, so there are infinitely many of these, giving the required contradiction. 
\end{proof}
\end{prel6}
\begin{prel7}\label{nestfin}
A mincut is nested with all but finitely many mincuts.
\begin{proof}
Suppose $K$ is a mincut and $L$ is a mincut not nested with $K$. By lemma \ref{ABthm} both $K,L$ are B-cuts, with components $C_1, C_2$ of $\Gamma - K$, $D_1, D_2$ of $\Gamma - L$. If $L\cap C_1$ were empty then by connectedness $C_1\subseteq D_1$ or $C_1\subseteq D_2$, both of which would imply that $K,L$ were nested. Similarly none of $C_2\cap L$, $D_1\cap K$, $D_2\cap K$ is empty. Then $L$ is a \xysep for some $x\in K\cap D_1$, $y\in K\cap D_2$. There are only finitely many such separators for each pair $x,y$ and only finitely many elements of $K$, so only finitely many such $L$ are possible. 
\end{proof}
\end{prel7}

%% file: crosscuts.tex
\begin{figure}[!ht]
\centering
\includegraphics[width=0.8\textwidth]{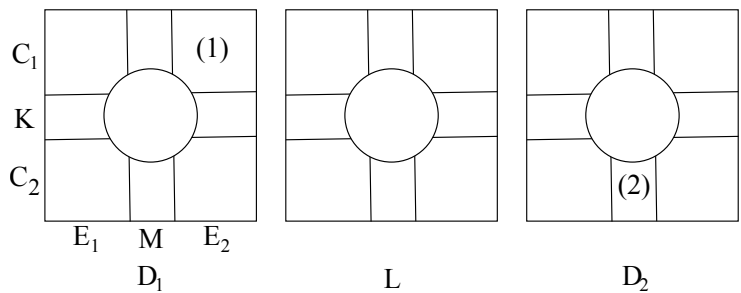}
\caption{The notation indicates that, for example, (1) is $C_1\cap D_1\cap E_2$ and (2) is $C_2\cap D_2\cap M$.}
\label{diagramexample}
\end{figure}

The complexity in the structure of mincuts comes from so-called `crossing' cuts, which we now define. 
\begin{crossdef}
Let $K, L$ be mincuts. Let $\mathcal{E}$ be the set of ends of $\Gamma$, and let $\mathcal{E}=K^{(1)}\sqcup K^{(2)}\sqcup\ldots\sqcup K^{(r)}$, $\mathcal{E}=L^{(1)}\sqcup L^{(2)}\sqcup\ldots\sqcup L^{(s)}$ be the partitions of $\mathcal{E}$ given by $K, L$ respectively. We say $[K], [L]$ \emph{cross} if, possibly after relabelling, $K^{(i)}\cap L^{(j)} \neq\emptyset$ for $i,j=1,2$. We write $K+L$.
\end{crossdef}
The following is a direct consequence of \ref{ABthm}, having removed slices from our graph.
\begin{Bcuts}
If $[K], [L]$ cross then $\Gamma - K$, $\Gamma - L$ have exactly two components.
\end{Bcuts}
Later we will show that crossing mincuts possess a cyclic structure. Initially however we shall just consider two or three crossing cuts.
\begin{cross1}
Let $[K], [L]$ be crossing classes of mincuts. Let $\Gamma - K = C_1\sqcup C_2, \Gamma - L = D_1\sqcup D_2$. Then $|C_1\cap L|=|C_2\cap L|=|D_1\cap K|=|D_2\cap K|$, i.e $K\cup L$ splits into four equal pieces, plus the `centre' $U = K\cap L$.
\begin{proof}
This follows from two applications of \ref{prel1}.
\end{proof}
\end{cross1}
In the case of edge cuts, one can also show that the centre $K\cap L$ is empty, but in the case of vertex cuts this fails to be true. As we will show in lemma \ref{magiclemma} below, the centre is in some sense distinguished, but this result must wait until we have placed some restrictions on the division of a graph produced by three cuts.\\
Let $K, L, M$ be mincuts with $K$ crossing $L$ and $L$ crossing $M$, and let $C_1\sqcup C_2$, $D_1\sqcup D_2$, $E_1\sqcup E_2$ be the components of $\Gamma-K$, $\Gamma - L$, $\Gamma - M$ respectively. A priori, these three cuts could divide $\Gamma$ into eight components each containing an end. The natural diagram with which to illustrate this would be a suitably divided cube. To produce this in 2D we divide the cube into three slices as shown in the following diagram, Figure \ref{diagramexample}.\\
We now rule out certain arrangements of ends of the graph.
\begin{cross2}\label{tetra}
It is not possible for each of $C_1\cap D_1\cap E_1, C_2\cap D_1\cap E_2, C_2\cap D_2\cap E_1, C_1\cap D_2\cap E_2$ to contain an end [or any arrangement obtained from this by relabellings].
\begin{proof}
Denote by $a,\ldots, u$ the cardinalities of the various subgraphs as shown below; the $\epsilon_i$ indicate the presence of ends.
\begin{figure}[htp]
\centering
\includegraphics[width=0.8\textwidth]{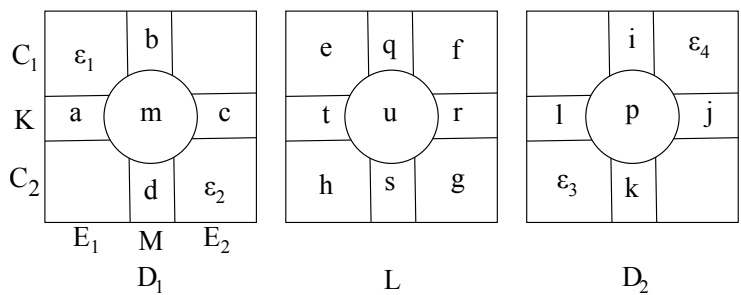}
\caption{Diagram for lemma \ref{tetra}}
\end{figure}
Let $n$ be the cardinality of a mincut. Then $|K| = |L| = |M| = n$, so:
\begin{IEEEeqnarray*}{rCl}
n&=&a+c+l+j+m+p+r+t+u\\
n&=&e+f+g+h+q+r+s+t+u\\
n&=&b+d+i+k+m+p+q+s+u
\end{IEEEeqnarray*}
We also have an end cut separating each $\epsilon_i$ from the others; this yields:
\begin{IEEEeqnarray*}{rCl}
n&\leq& a+b+e+m+q+t+u\\
n&\leq& c+d+g+m+r+s+u\\
n&\leq& k+l+h+p+s+t+u\\
n&\leq& i+j+f+p+q+r+u
\end{IEEEeqnarray*}
Sum these four:
\begin{IEEEeqnarray*}{rCl}
4n&\leq& (a+c+l+j+m+p+r+t+u)\\
&& + (e+f+g+h+q+r+s+t+u)\\
&& + (b+d+i+k+m+p+q+s+u) + u\\
&=& 3n + u
\end{IEEEeqnarray*}
whence $u = n$, everything else vanishes, and $K=L=M$, each separating the graph into at least four components, contradicting lemma \ref{tetra}.
\end{proof}
\end{cross2}
Note that this result implies that the three cuts split the graph into at most six components containing ends. Since $K$ crosses $L$ there are at least four such components. A quick exercise in filling in corners with ends subject to the crossings and lemma \ref{tetra} shows that, following relabellings, $C_1\cap D_1\cap E_1, C_2\cap D_1\cap E_2, C_1\cap D_2\cap E_1, C_2\cap D_2\cap E_2$ all contain ends [with possibly other corners also]. We can now prove:
\begin{cross3}\label{magiclemma}
Let $K, L, M$ be mincuts with $K$ crossing $L$ and $L$ crossing $M$ (in particular, if $K$ is equivalent to $M$). Then $K\cap L = L\cap M$, and $L\cap C_1 = L\cap E_1, L\cap C_2 = L\cap E_2$.
\begin{proof}
Retain the notations of the previous lemma. Again since $K, L, M$ are mincuts, 
\begin{IEEEeqnarray*}{rCl}
n&=&a+c+l+j+m+p+r+t+u\\
n&=&e+f+g+h+q+r+s+t+u\\
n&=&b+d+i+k+m+p+q+s+u
\end{IEEEeqnarray*}
and again, considering end cuts separating a corner containing an end $\epsilon_i$ from the others, we have:
\begin{IEEEeqnarray*}{rCl}
n&\leq& a+b+e+m+q+t+u\\
n&\leq& c+d+g+m+r+s+u\\
n&\leq& i+l+e+p+q+t+u\\
n&\leq& j+k+g+p+r+s+u
\end{IEEEeqnarray*}
Summing these,
\begin{IEEEeqnarray*}{rCl}
4n&\leq& (a+c+l+j+m+p+r+t+u)\\
&& + (b+d+i+k+m+p+q+s+u)\\
&& + 2e + 2g +q+ r+s+t + 2u\\
&=& 2(e+f+g+h+q+r+s+t+u)\\
&& + 2n - 2f-2h-q-r-s-t\\
&=& 4n - (2f + 2h +q+r+s+t)
\end{IEEEeqnarray*}
\begin{figure}[ht]
\centering
\includegraphics[width=0.8\textwidth]{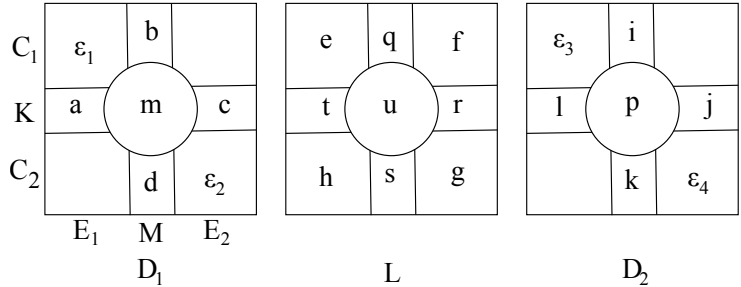}
\caption{Diagram for lemma \ref{magiclemma}}
\end{figure}
whence $f=h=q=r=s=t=0$, so that $K\cap L = L\cap M$.
\end{proof}
\end{cross3}
\begin{cross4}
A cut is crossed by at most finitely many cuts. 
\begin{proof}
If two cuts cross they are not nested, so this follows directly from lemma \ref{nestfin}
\end{proof}
\end{cross4}

%% file: halfcuts.tex
It follows from the last section's results that each mincut in a crossing system can be decomposed into three pieces; two `half-cuts' and a `centre'. We now prove some facts about these half-cuts, which enable us to arrange the half-cuts on a circle.
\begin{halfdefn}
If $K,M$ are mincuts (more properly, classes of mincuts under $\sim$, but we will often pass over this technicality), we write $K$\#$L$ if there are mincuts $K=L_0, L_1, \ldots, L_n = M$ such that $L_0 + L_1 +\ldots + L_n$; that is, $L_0$ crosses $L_1$, $L_1$ crosses $L_2$ and so on. $L_0$ may or may not cross $L_2$. \# is an equivalence relation on $\sim$-classes of mincuts, decomposing these into equivalence classes, which we call \#-classes.
\end{halfdefn}
By lemma \ref{magiclemma}, elements $K$ of a \#-class have a unique decomposition $K= K_1 \cup U \cup K_2$ where if $K+L$ then $K\cap L = U$ and $K_1, K_2$ are in different components of $\Gamma - L$. From the same lemma, this $U$ is the same for all cuts in the \#-class, we call it the \emph{centre} of the \#-class. Also, $|K_1|=|K_2|$ and this cardinality is again the same across the class. The $K_i$ are called the \emph{half-cuts} of the \#-class. We now prove a series of lemmas clarifying the structure of a \#-class and its half-cuts.
\begin{half1}\label{lhalf1}
If $K,M$ are mincuts in the same \hashc then either $K+M$ or there is an $L$ in this \hashc such that $K+L+M$; that is, $K+L$ and $L+M$.
\begin{proof}
By definition we have a sequence of cuts $K=L_0, L_1, \ldots, L_n = M$ such that $L_0 + L_1 +\ldots + L_n$. Take a shortest such sequence, and suppose $n\geq 3$. We will show we can find a shorter sequence. Without loss of generality we can assume $n=3$. Let $\mathcal{E} = L_i^{(1)}\cup L_i^{(2)}$ be the partition induced by $L_i$. The fact that $K$ does not cross $L_2$, and similar facts, give us after relabelling:
\begin{IEEEeqnarray*}{rCllrCl}
K^{(1)}&\subseteq & L_2^{(1)}&,\qquad&L_2^{(2)}&\subseteq & K^{(2)}\\
M^{(2)}&\subseteq & L_1^{(2)}&,\qquad&L_1^{(1)}&\subseteq & M^{(1)}\\
\end{IEEEeqnarray*} 
whence the crossings give us that each of $M^{(2)}\cap K^{(2)}, M^{(1)}\cap K^{(2)}, M^{(1)}\cap K^{(1)}$ is non-empty. Hence $K+M$ unless $M^{(2)}\cap K^{(1)}$ is empty, hence $K^{(1)}\subseteq M^{(1)}$. It is this that allows us to place the ends $\epsilon_3, \epsilon_6$ in Figure \ref{figure31}, and hence to conclude that $K+ (L_{11}\cup U\cup L_{22}) + M$, where for instance $L_{11}$ is the half-cut of $L_1$ lying in the $K^{(1)}$-component of $\Gamma - K$.
\begin{figure}[htp]
\centering
\includegraphics[width=0.8\textwidth]{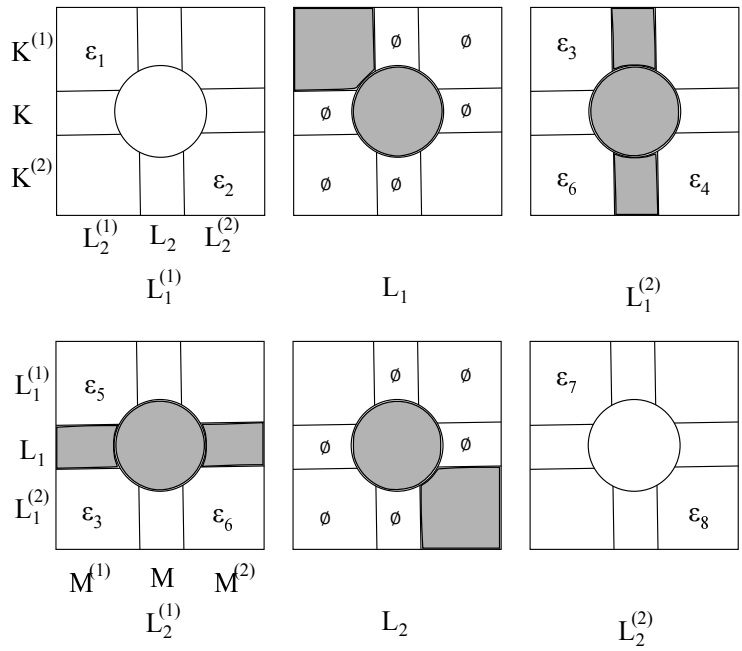}
\caption{Diagram for lemma \ref{lhalf1}. The cut $(L_{11}\cup U\cup L_{22})$ is shown shaded in both diagrams.}
\label{figure31}
\end{figure}
\end{proof}
\end{half1}
\begin{half2}
There are only finitely many cuts in a \hashc and hence only finitely many half-cuts.
\end{half2}
\begin{half3}\label{lhalf3}
Let $K_1, M_1$ be half-cuts in the same \#-class. Then $K_1\cup U\cup M_1$ is a mincut iff there are mincuts $K', M'$ containing $K_1, M_1$ as half-cuts respectively such that $K'+M'$.
\begin{proof}
One direction is obvious. For the other, pick $K_2, M_2$ such that $K=K_1\cup U\cup K_2$, $M=M_1\cup U\cup M_2$ are cuts in this \#-class. Then either $K+M$, in which case we are done, or there is $L$ such that $K+L+M$. $K_1\cup U\cup M_1$ is a cut, hence we have $\epsilon_5$ in Figure \ref{figure34}. Then $K_1\cup U\cup M_2$, $K_2\cup U\cup M_1$ cross.
\begin{figure}[htp]
\centering
\includegraphics[width=0.8\textwidth]{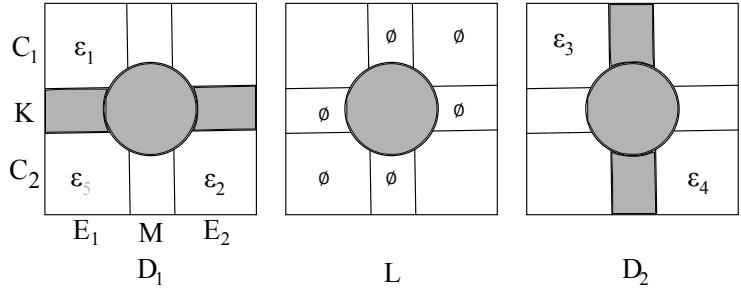}
\caption{Diagram for lemma \ref{lhalf3}. The cut $(K_1\cup U\cup M_2)$ is shown shaded.}
\label{figure34}
\end{figure}
\end{proof}
\end{half3}
\begin{halfdef2}
Two half-cuts $K_1, L_1$ in the same \hashc are \emph{equivalent} if whenever $K_2$ is a half-cut such that $K_1\cup U \cup K_2$ is a mincut, then $L_1 \cup U \cup K_2$ is an equivalent cut and vice versa.\\
Two half-cuts $K_1, L_1$ in the same \hashc are \emph{\qe} if there is a half-cut $K_2$ such that $K_1\cup U \cup K_2$ is a mincut and $L_1 \cup U \cup K_2$ is an equivalent cut. 
\end{halfdef2}
\begin{half4}\label{half4}
If two half-cuts $K_1, M_1$ form a cut then they are not \qe.
\begin{proof}
Let $K=K_1\cup U \cup M_1$ be the cut formed by hypothesis. Let $L_1$ be some other half-cut; we will show that $K_1\cup U\cup L_1$ is not equivalent to $L_1\cup U\cup M_1$ as cuts, hence that $K_1, M_1$ are not quasi-equivalent. Let $L_2$ be a half-cut such that $L= L_1\cup U \cup L_2$ is in the \#-class. If $L+K$ then the result is clear. If not, there is a mincut $N$ such that $K+N+L$; without loss of generality take $K_1, L_1$ to be in the same component of $\Gamma - N$. Then $L_1\cup U\cup M_1$ is a cut, and from the diagram we see that either $K_1\cup U \cup L_1$ is not an end cut or it is not equivalent to $L_1\cup U \cup M_1$.
\begin{figure}[htp]
\centering
\includegraphics[width=0.8\textwidth]{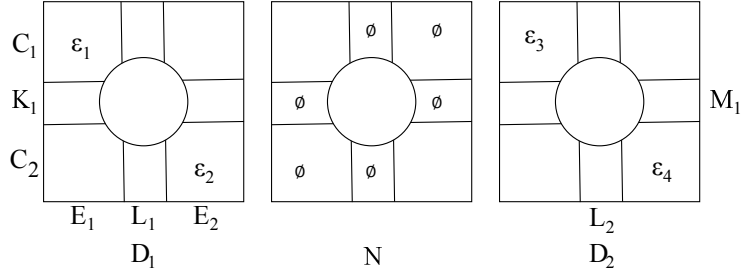}
\label{figure36}
\caption{Diagram for lemma \ref{half4}}
\end{figure}
\end{proof}
\end{half4}
\begin{half5}\label{half5}
Let $K=K_1\cup U\cup K_2$ be a cut in the \hashc and let $M_1$ be a half-cut in the same class not \qe\ to either $K_1, K_2$. Then there is $M_2$ such that $M_1\cup U\cup M_2$ is a cut crossing $K$. Hence $K_1\cup U\cup M_1$ and $K_2\cup U\cup M_2$ are cuts.
\begin{proof}
Let $M_2$ be a half-cut such that $M=M_1\cup U\cup M_2$ is a cut of the class. If $K+M$ we are done. Otherwise, there is a cut $L$ with $K+L+M$. After possibly relabelling the $K_i$, we can assume that $K_1, M_1$ are in the same component of $\Gamma - L$. If there is an end in $C_1\cap D_1\cap E_2$ then $M_1\cup U\cup L_2$ is a cut crossing $K$. If there is an end in $C_2\cap D_1\cap E_1$ then $M_1\cup U\cup L_1$ is a cut crossing $K$. If neither of these happens, then $K_2\cup U\cup M_1$ is equivalent to $K_2\cup U\cup K_1$, a contradiction (we note that these cuts are genuinely equivalent, since the presence of `links' such as $L_1$ guarantees that ends which appear to be connected up really are).
\begin{figure}[htp]
\centering
\includegraphics[width=0.8\textwidth]{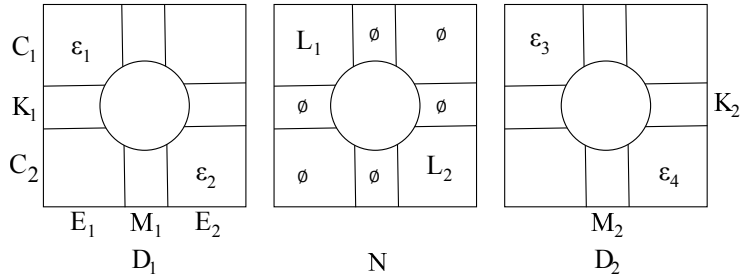}
\caption{Diagram for lemma \ref{half5}}
\label{figure37}
\end{figure}
\end{proof}
\end{half5}
\begin{half6}\label{welldefine}
Quasi-equivalence is an equivalence relation. If $K_1, L_1$ are \qe\ and $L_2$ is a half-cut such that $L_1\cup U\cup L_2$ is in the \hashc then $K_1\cup U\cup L_2 \sim L_1\cup U\cup L_2$.
\begin{proof}
Let $K_2$ be such that $K=K_1\cup U\cup K_2$ is in the \#-class. The cut $K$ does not cross $L=L_1\cup U\cup L_2$ since in this case $K_1\cup U\cup L_1$ would be a cut, so $K_1, L_1$ are not \qe\ by lemma \ref{half4}, giving a contradiction. Then there is an $N$ such that $K+N+L$. Again, $K_1\cup U\cup L_1$ is not a cut, so there are no ends in certain corners as indicated. Then $K_1\cup U\cup L_2\sim L_1\cup U\cup L_2$ as required, noting again that ends which appear connected actually are so that the cuts are genuinely equivalent.\\
As for quasi-equivalence being an equivalence relation, it is clearly symmetric and reflexive. If $M_1$ is another half-cut \qe\ to $L_1$, then by the above
\begin{equation*}
K_1\cup U\cup L_2\sim L_1\cup U\cup L_2\sim M_1\cup U\cup L_2
\end{equation*}
so $K_1, M_1$ are \qe.
\begin{figure}[htp]
\centering
\includegraphics[width=0.8\textwidth]{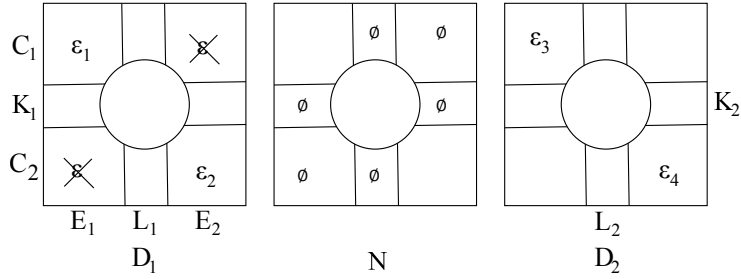}
\caption{Diagram for lemma \ref{welldefine}}
\label{figure38}
\end{figure}
\end{proof}
\end{half6}
\begin{half7}\label{issep}
Let $K= K_1\cup U\cup K_2$ be a cut in the \hashc and let $L_1, M_1$ be half-cuts in the same class not \qe\ to $K_1, K_2$. Then either $L_1 \cup U \cup M_1$ is a cut crossing $K$ or $L_1, M_1$ are contained in the same component of $\Gamma - K$.
\begin{proof}
By lemma \ref{half5} we can complete $L_1$ to a cut crossing $K$, so that $L_1$ separates some ends of a component of $\Gamma - K$; similarly for $M_1$. If the two half-cuts are in different components, then $L_1\cup U\cup M_1$ is a cut crossing $K$ provided it separates $\Gamma - K$ into two components, as indeed it must.
\end{proof}
\end{half7}

%% file: separations.tex
We now turn our attention to demonstrating that the half-cuts of a system have a cyclic structure. We will do this by showing that they satisfy a certain axiomatic system, which implies that they can be arranged cyclically in a fashion compatible with their cut structure. This axiomatic structure is taken from a 1935 paper of Huntington \cite{huntington1935inter}.
\begin{sepdef}
A separation relation on a set $Z$ is a relation $R\subseteq Z^4$, satisfying the following axioms. We write $abcd$ if $(a,b,c,d)\in R$.
\begin{enumerate}
\item If $abcd$ then $a, b, c, d$ are distinct.
\item There are $a, b, c, d$ such that $abcd$, i.e. $R\neq\emptyset$.
\item If $abcd$ then $bcda$
\item If $abcd$ then $\neg(abdc)$
\item There are $a, b, c, d$ such that $abcd$ and $dcba$.
\item If $abcd$ and $x\in Z$ is another element then either $axcd$ or $abcx$.
\end{enumerate}
\end{sepdef}
\begin{sep1}\label{seplemmas}
Let $Z$ be a set equipped with a separation relation. Then:
\begin{enumerate}
\item If $a,b,c,d\in Z$ are distinct, then at least one of the twenty-four tetrads $abcd, abdc,\ldots,dcba$ is true. 
\item If $abcd$ then $dcba$.
\item If $abxc$ and $abcy$ then $abxy$.
\item If $abcx$ and $abcy$ then $abxy$ or $abyx$.
\item If $abcx$ and $abcy$ then $acxy$ or $acyx$.
\end{enumerate}
where in the last three statements distinct letters are assumed to represent different elements of $Z$.
\end{sep1}
Proofs can be found in \cite{huntington1932postulates} along with further similar propositions.
\begin{sep2}
Let $Z$ be a finite set with a separation relation. For each $z$ there are unique $a,b$ such that for all $c\in  Z - \{z,a,b\}$, $azbc$. We call these the elements adjacent to $z$.
\begin{proof}
We approach existence by induction. For $|Z|=4$ the result is trivial. Assume it is true for all separation relations with $|Z|=n$, and suppose $|Z|=n+1$. Remove an element $d$ of $Z$ not equal to $z$ to leave a smaller separation relation, and let $a,b$ be the elements adjacent to $z$ in this new relation, so that for all $c\in Z-\{z,a,b,d\}$, $azbc$.\\
By lemma \ref{seplemmas}, one of $azbd$, $adzb$, $azdb$ holds. If $azbd$ holds then $a,b$ are adjacent to $z$ in $Z$. If not, without loss of generality, $azdb$. We claim $a,d$ are adjacent to $z$ in $Z$. For by lemma \ref{seplemmas} above, if $c\in Z-\{z,a,d,b\}$ then $azdb$ and $azbc$ imply $azdc$.\\
For the uniqueness part, suppose there are two such pairs $a_1, b_1$, $a_2, b_2$. If any of these coincide we have an immediate contradiction to part 4 of the definition of the relation. So suppose they are all distinct. Then $a_1 z b_1 a_2, a_1 z b_1 b_2$ imply $a_1 z a_2 b_2$ or $a_1 z b_2 a_2$, both of which contradict $a_2 z b_2 a_1$.
\end{proof}
\end{sep2}
\begin{sep3}
Let $Z$ be a finite set equipped with a separation relation. Then there is a map $F:Z\to S^1$ such that for $a, b,c,d\in Z$, $abcd \iff F(b), F(d)$ lie in different components of $S^1 - \{F(a), F(c)\}$, i.e. $Z$ is isomorphic to a finite subset of the circle under its natural separation relation.
\begin{proof}
By induction. Pick an element $z\in Z$ and take a separation-preserving map $\bar{F}:Z-\{z\}\to S^1$. By the previous lemma, there are elements $a,b$ of $Z$ adjacent to $z$. We will map $z$ to the circle by placing it between $\bar{F}(a), \bar{F}(b)$, but first we must show these are adjacent on the circle. If not, there are $c,d$ so that $\bar{F}(a)\bar{F}(c)\bar{F}(b)\bar{F}(d)$, whence $acbd$. But $azbc, azbd$ imply $abcd$ or $abdc$, both contradicting $acbd$. So $\bar{F}(a), \bar{F}(b)$ are adjacent on the circle, and we can define $F:Z\to S^1$ by setting $F=\bar{F}$ on $Z-\{z\}$ and $F(z)$ to lie between $\bar{F}(a), \bar{F}(b)$ on the circle.\\
A full proof that this $F$ works would be lengthy and uninformative, so we just indicate the main steps; the remainder is just use of axioms and \ref{seplemmas}. We inherit from $\bar{F}$ that any relations not involving $z$ are preserved. Let $zABC$ be another relation and suppose $A,B,C$ are distinct from $a,b$; the other cases are easier. Then we have $azbA, azbB, azbC, zABC$ from which we deduce $aABC, bABC$. These relations carry over to the circle under $F$, as do $azbA, azbB, azbC$ by construction. From these relations on the circle we then find $F(z)F(A)F(B)F(C)$.
\end{proof}
\end{sep3}
\begin{sepdef2}
Let $Z$ be the set of quasi-equivalence classes of half-cuts of a \#-class. We define a separation relation $R$ on $Z$ by setting $(a,b,c,d)\in R$ if and only if $ac + bd$, where $ac$ denotes the cut $K_1\cup U\cup L_1$ and $K_1, L_1$ are representatives of $a, c$ and so on.
\end{sepdef2}
\begin{sep4}
This is well-defined, i.e. it does not matter which representatives of quasi-equivalence classes we choose. Furthermore, it is a bona fide separation relation.
\begin{proof}
Well-definedness follows immediately from lemma \ref{welldefine}. Parts 1-5 of the definition of a separation relation are trivial. For part 6, by lemma \ref{issep} either $abcx$ or $b,x$ are in the same component of $\Gamma - {ac}$; and either $axcd$ or $x,d$ are in the same component of $\Gamma - ac$. But $b,d$ are in different components of $\Gamma- ac$ so one of $abcx, axcd$ holds.
\end{proof}
\end{sep4}
Hence we have:
\begin{sep5}\label{cyclic}
To each \hashc we can associate a cycle where each vertex represents a quasi-equivalence class and each cut of the \hashc is associated to a vertex cut of the cycle, with the notions of crossing preserved.
\end{sep5}

%% file: classstructure.tex
We are now in a position to characterize the structure of a \#-class. Let $[K_1]_q$ denote the quasi-equivalence class of a half-cut $K_1$. From lemma \ref{cyclic} there are two quasi-equivalence classes adjacent to $[K_1]_q$ in this \#-class. If $L_1$ is a half-cut in the \hashc not in $[K_1]_q$ or either quasi-equivalence class adjacent to it, then lemma \ref{welldefine} implies that 
\begin{equation*}
K_1\cup U\cup L_1 \sim K'_1\cup U\cup L_1
\end{equation*}
for all $K'_1\in [K_1]_q$. So only the two quasi-equivalence classes adjacent to $[K_1]_q$ can contain $L_1$ such that $K=K_1\cup U\cup L_1$ and $K'=K'_1\cup U\cup L_1$ are not equivalent for $K'_1\in [K_1]_q$.\\
How can these cuts be non-equivalent? We recall that by minimality every component left by a mincut is connected to every element of that cut. Thus in the ``larger component" left by the cut, i.e. the one containing half-cuts in the same class, every vertex is connected to the half-cuts in this ``component", which is thus genuinely connected. Thus one part of the partitions induced by $K, K'$ is the same. The others can only differ if at least one of the cuts splits $\Gamma$ into more than two parts, hence the ``smaller" component into more than one part. Suppose $K$ intersects one of the ``smaller" components of $\Gamma - K'$. Then each end not in the larger component of $\Gamma - K'$ is connected to each vertex of the part of $K$ in the ``smaller" component, hence $K'$ splits $\Gamma$ into exactly two components. If $K$ does not intersect one of the ``smaller" components of $\Gamma - K'$, then since $K\neq K'$ and $K, K'$ have the same cardinality, $K'$ intersects one of the ``smaller" components of $\Gamma - K$, whence $K$ splits $\Gamma$ into exactly two components.\\
Hence, having chosen $L_1$ there are at most two equivalence classes of cuts formed from $L_1$ and $[K_1]_q$. By symmetry, there are at most two equivalence classes of cuts formed from $K_1$ and $[L_1]_q$. From these discussions it follows that for each quasi-equivalence class adjacent to $[K_1]_q$ there are at most two equivalence classes of cuts formed by these two classes; one producing a split of $\Gamma$ into two components, the other more. Hence there are at most four equivalence classes of half-cuts within $[K_1]_q$.\\
We now define the structure by which we model the \#-class. For edge cuts this would be a simple cycle, but here we need extra complexity to deal with the possibility of splitting the graph into more than two components.
\begin{ringdef}\label{rindef}
A \emph{ring} is constructed as follows. Take a finite cycle of vertices and attach to each edge some number of triangles by identifying an edge of the triangle with the edge of the cycle. The vertex of the triangle not included in the original cycle is called an \emph{anchor}.\\
\end{ringdef}
\begin{nvertexdef}
An \emph{$n$-vertex} will be a copy of the complete graph on $n$ vertices; we will say it is connected to a vertex if there is an edge from the vertex to each constituent vertex of the $n$-vertex. We will depict a 3-vertex as a triangle and only draw one edge from it to each vertex to which it is connected. 
\begin{figure}[htp]
\centering
\includegraphics[width=0.8\textwidth]{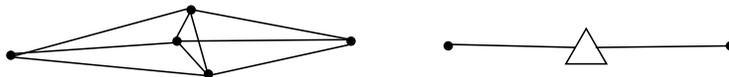}
\caption{A 3-vertex connected to two 1-vertices, and the schematic representation of this.}
\label{figureanchor}
\end{figure}
\begin{figure}[htp]
\centering
\includegraphics[width=0.4\textwidth]{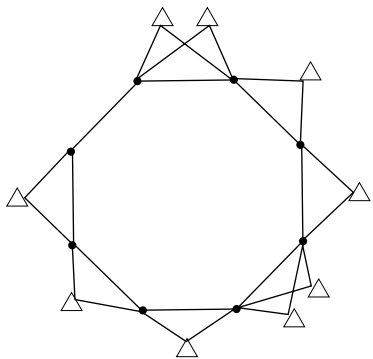}
\label{figurering}
\caption{a ring, with the anchors replaced by 3-vertices.}
\end{figure}
\end{nvertexdef}
We now associate to each \hashc an appropriate ring encoding the cuts formed by half-cuts in the class. First use proposition \ref{cyclic} to form a cycle with one vertex for each quasi-equivalence class. For each pair of adjacent quasi-equivalence classes, find half-cuts in those classes separating $\Gamma$ into as many components as possible, and attach one fewer anchors than this between the two classes in the cycle (one fewer to account for the ``large" component). If a quasi-equivalence class contains more than one equivalence class, insert an extra vertex into the cycle here. If we ``thicken'' up the anchors to 3-vertices to remove cut-points, there is now a bijective correspondence between equivalence classes of cuts formed from half-cuts of the \hashc and equivalence classes of cuts of the ring, where we treat the anchors as ends for the purpose of equivalence etc.
 \begin{figure}[ht]
\centering
\includegraphics[width=0.8\textwidth]{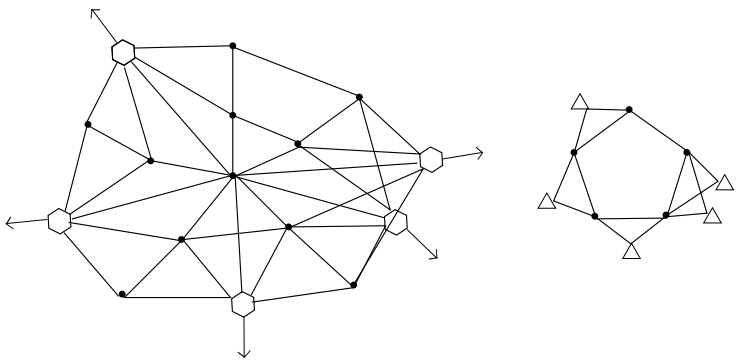}
\label{figureringexample}
\caption{a \hashc and its associated ring. Hexagons represent 6-vertices and arrows ends.}
\end{figure}

%% file: pretrees.tex
We now proceed towards the central theorem of the paper. First we seek to impose a tree structure on the \#-classes and the other cuts, then we will re-introduce the extra complexity. We will do this using pretrees, which we now define.
\begin{pretdef1}
Let $\mathcal{P}$ be a set and let $R\subseteq \mathcal{P}^3$ be a ternary relation on $\mathcal{P}$. If $(x,y,z)\in R$ then we write $xyz$ and say $y$ is between $x,z$. A set $\mathcal{P}$ equipped with this relation is a \emph{pretree} if the following hold:
\begin{enumerate}
\item If $xyz$ then $y\neq x,z$, and there are no $x,y$ such that $xyx$.
\item If $xyz$ then $zyx$.
\item For all $x,y,z$, if $xyz$ then $\neg(xzy)$.
\item If $xzy$ and $w\neq x,y,z$ then $xzw$ or $yzw$.
\end{enumerate}
If there is no $z$ such that $xzy$ we say, $x,y$ are adjacent.\\
A pretree is called \emph{discrete} if for any $x,y\in\mathcal{P}$ there are at most finitely many $z$ such that $xzy$.
\end{pretdef1}
It should perhaps be noted that despite us using the word `between' this is not a betweenness relation in the usual sense of the word as for example in \cite{huntington1935inter}. Let $\mathcal{P}$ be a discrete pretree. We will describe briefly how to pass from $\mathcal{P}$ to a tree; a fuller description may be found in \cite{bowditch1999treelike}.\\
We call a subset $H$ of $\mathcal{P}$ a \emph{star} if all $a,b\in H$ are adjacent. We now define a tree $T$ as follows:
\begin{equation*}
V(T) = \mathcal{P}\cup \{\text{maximal stars of } \mathcal{P}\}
\end{equation*}
\begin{equation*}
E(T) = \{(v,H): v \in\mathcal{P}, v\in H, H \text{ a maximal star}\}
\end{equation*}
We show that $T$ is indeed a tree. If $x,y\in\mathcal{P}$ then by discreteness there are only finitely many $z$ between $x,y$. From among these $z$ we can then find $z_1,\ldots,z_n$ such that $x$ is adjacent to $z_1$, $z_i$ is adjacent to $z_{i+1}$ and $z_n$ is adjacent to $y$, giving a path in $T$ from $x$ to $y$. Hence $T$ is connected.\\
If $T$ contains a circuit then there are $x_1,\ldots,x_n$ in $\mathcal{P}$ such that $x_i$ is adjacent to $x_{i+1}$ but not to $x_{i+2}$ for each $i\in\mathbb{Z}_n$. Then there is $y$ such that $x_i y x_{i+2}$. If $y\neq x_{n+1}$ then either $x_i y x_{i+1}$ or $x_{i+1} y x_{i+2}$ both of which are forbidden. So $x_i x_{i+1} x_{i+2}$. We claim $x_1 x_i x_{i+1}$ holds for all $i\leq n$ by induction. Since $x_1 x_{i-1} x_{i}$ and $x_{i-1}\neq x_{i+1}$, $x_1 x_{i-1} x_{i+1}$ holds or we have a contradiction. Since $x_{i-1} x_i x_{i+1}$ and $x_1 \neq x_i$, either $x_1 x_i x_{i+1}$ or $x_{i-1} x_i x_{1}$; so to avoid contradiction, $x_1 x_i x_{i+1}$. But then we have $x_1 x_{n-1} x_n$; but $x_1, x_n$ were supposed to be adjacent. The contradiction means $T$ is a tree.\\
We now prove some lemmas which will allow us to define a pretree of cut classes.
\begin{pretdef2}
We call a mincut \emph{isolated} if it does not cross any mincut, hence is not contained in any \#-class.\\
A cut is a \emph{corner cut} of a \hashc if it is (equivalent to) a cut formed from two half-cuts of the class but is not itself in the class. We call a mincut \emph{totally isolated} if it does not cross any mincut and is not a corner cut of any \#-class.
\end{pretdef2}
\begin{pret1}
Corner cuts are isolated.
\begin{proof}\label{pret1}
Let $Q$ be a \hashc and let $K=K_1\cup U\cup K_2$ be a corner cut of $Q$. Suppose there is a cut $L$ with $K+L$. Then $L$ separates some ends of each component  of $\Gamma - K$. Let $M_1, M_2$ be half-cuts in $Q$ adjacent to $K$, with $K_1$ adjacent to $M_2$ and $K_2$ adjacent to $M_1$, with no quasi-equivalences present. Either $L$ crosses $K_1\cup U\cup M_1$ or all ends of the component of $\Gamma-K_1\cup U\cup M_1$ containing $M_2$ are in the same component of $\Gamma - L$, whence $L$ crosses $K_2\cup U\cup M_2$. So $L$, hence $K$, are in the \hashc $Q$. Contradiction.
\begin{figure}[htp]
\centering
\includegraphics[width=0.4\textwidth]{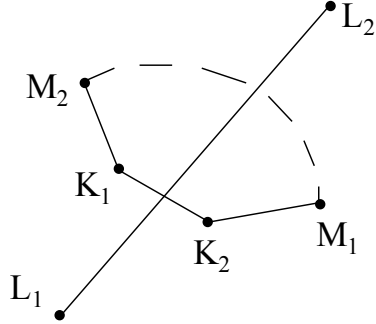}
\caption{Diagram for lemma \ref{pret1}.}
\end{figure}
\end{proof}
\end{pret1}
Each \hashc $Q$ induces two partitions 
\begin{eqnarray*}
\mathcal{E} = Q^{(1)}\sqcup\ldots\sqcup Q^{(m)}\\
\mathcal{E} = \bar{Q}^{(1)}\sqcup\ldots\sqcup \bar{Q}^{(m')}
\end{eqnarray*}
of the ends of $\Gamma$. In one partition, which we call the fine partition and denote without bars, each member of the partition corresonds to one of the anchors in the ring representing $Q$; and for each $Q^{(i)}$ there is a corner cut of $Q$ separating $Q^{(i)}$ from all the other $Q^{(j)}$. For the other partition, the coarse partition, we identify those \pn{Q}{i}\ together which lie between the same two adjacent half-cuts. Then in the coarse partition we can distinguish between members \pn{\bar{Q}}{i}\ using only cuts properly in the \hashc $Q$; for the fine partition we may need corner cuts also. We recall also that a cut $K$ also gives a partition of the ends of $\Gamma$:
\begin{equation*}
\mathcal{E} = K^{(1)}\sqcup\ldots\sqcup K^{(n)}
\end{equation*}
\begin{pret2}\label{ticc}
Given a cut $K$ and a \hashc $Q$, with $K$ neither in $Q$ nor a corner cut of it, there are $i,j$ such that all $Q^{(k)}$ except $Q^{(i)}$ are contained in $K^{(j)}$, i.e.
\begin{equation*}
\coprod_{k\neq i} Q^{(k)} \subseteq K^{(j)}
\end{equation*}
or:
\begin{equation*}
\coprod_{k\neq j} K^{(k)} \subseteq Q^{(i)}
\end{equation*}
We say $K$ \emph{divides} $Q^{(i)}$.
\begin{proof}
Suppose $K$ is an A-cut. Then it is nested with every cut and corner cut of $Q$, hence the result.\\
Otherwise $K$ is a B-cut, separating $\Gamma$ into two components. If the result is not true, then both $K^{(i)}$ intersect at least two $Q^{(i)}$. \\
Suppose a $Q^{(i)}$ intersects both $K^{(i)}$. Let $M$ be the corner cut of $Q$ splitting off \pn{Q}{i}. If $M$ is a B-cut, then $K+M$ giving a contradiction. Otherwise, $M$ is nested with $K$, whence a \pn{K}{i} is contained in a \pn{Q}{i}, again giving a contradiction.\\
If both $K^{(i)}$ contain two $Q^{(i)}$ not between two adjacent half-cuts (in the ring representing $Q$) we can find a cut of $Q$ crossing $K$. So for say $K^{(1)}$ all the $Q^{(i)}$ contained in $K^{(1)}$ lie between two adjacent half-cuts of $Q$. Let $M$ be the corner cut corresponding to these half-cuts.\\
$M$ is necessarily an A-cut, hence is nested with $K$. As in the discussion of quasi-equivalent cuts earlier, $K$ can only intersect the ``large" component of $\Gamma-M$, that containing the other half-cuts of $Q$; conversely $M$ does not intersect the ``large" component $C_1$ of $\Gamma - K$. Pick another half-cut $L_1$ in $Q$ with $L= M_1\cup U\cup L_1$ a cut of $Q$. The cut $L$ also only intersects the ``large" component $E_1$ of $\Gamma - M$. With suitable labelling of the $L^{(i)}$, we have:
\begin{IEEEeqnarray*}{rCl}
L^{(1)}&\subseteq &M^{(1)}\\
L^{(1)}&\subseteq &K^{(1)}\\
K^{(2)}&\subseteq &L^{(2)}\\
K^{(1)}&=&M^{(1)}\\
\end{IEEEeqnarray*}
Hence we have the arrangement shown in Figure \ref{figure63}.
\begin{figure}[htp]
\centering
\includegraphics[width=0.8\textwidth]{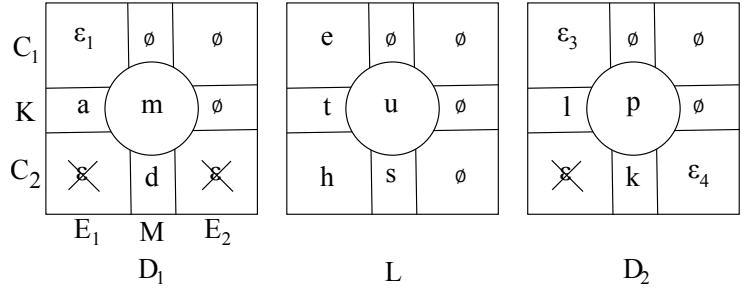}
\caption{Diagram for lemma \ref{ticc}.}
\label{figure63}
\end{figure}
Then we have 
\begin{IEEEeqnarray*}{rCl}
a+m+p+l+t+u&=&n\\
d+m+p+k+s+u&=&n\\
e+t+u+h+s&=&n\\
a+m+e+t+u&\geq& n\\
p+l+e+t+u&\geq& n\\
p+k+u+s&\geq& n
\end{IEEEeqnarray*}
where $n$ is the cardinality of a mincut. Immediately $d=m=0$. Furthermore, since $L_2, M_1, M_2$ are half-cuts in the same class, $s=e+t+h$. Then
\begin{IEEEeqnarray*}{rCl}
2n &\leq& (a+e+t+u)+(p+l+e+t+u)\\
&=& (a+p+l+t+u) + (e+e+t+u)\\
&\leq& n + e+t+h+s+u\\
&=&2n
\end{IEEEeqnarray*}
Then all the inequalities are equalities, hence $t=h=0$, and $K$ decomposes into $U$ together with two equal half-cuts; and choosing $L_1$ appropriately we find that these half-cuts are quasi-equivalent to the half-cuts of $M$, so $K$ was a corner cut of $Q$.
\end{proof}
\end{pret2}
\begin{pret3}\label{cccc}
Given two \#-classes $Q,R$, all cuts in $R$ divide the same $Q^{(i)}$
\begin{proof}
Note first that the cuts in $R$ do divide a $Q^{(i)}$ because they are not isolated, hence not a corner cut, and are not in $Q$. Suppose $K\in R$ divides \pn{Q}{i}\ and $L\in R$ divides \pn{Q}{j}, with $i\neq j$.\\
If there is one, take a cut $M$ crossing $K$ and $L$. We have \pn{K}{2}$\subseteq$\pn{Q}{i} and \pn{L}{2}$\subseteq$\pn{Q}{j} so $M$ contradicts lemma \ref{ticc}.\\
Then $K+L$. Take a cut $M\in Q$ separating \pn{Q}{i}\ from \pn{Q}{j}, and let $N=K_2\cup U\cup L_2$. The cut$N$ separates some ends of \pn{Q}{i}\ and of \pn{Q}{j}; it is not nested with $M$ hence is a B-cut and crosses $M$ giving a contradiction.
\begin{figure}[htp]
\centering
\includegraphics[width=0.3\textwidth]{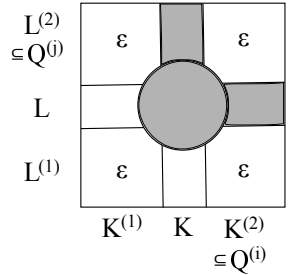}
\caption{Diagram for lemma \ref{cccc}. The cut $N$ is shown shaded.}
\label{figure75}
\end{figure}
\end{proof}
\end{pret3}
\begin{pret4}\label{titi}
Given two totally isolated cuts $K, L$, $K$ divides only one \pn{L}{i}, i.e. there are $i,j$ such that
\begin{IEEEeqnarray*}{rCl}
\coprod_{k\neq i} \pnm{L}{k}&\subseteq &\pnm{K}{j}\\
\coprod_{k\neq j} \pnm{K}{k}&\subseteq &\pnm{L}{i}\\
\end{IEEEeqnarray*}
\begin{proof}
If $K,L$ are nested, the result is immediate. If not, they are both B-cuts, when the result follows since they do not cross.
\end{proof}
\end{pret4}
We now define a pretree encoding the mincuts of $\Gamma$. Let $\mathcal{P}$ be the set of all \#-classes of $\Gamma$ and all equivalence classes of totally isolated cuts of $\Gamma$. Given $x,y,z\in\mathcal{P}$, we say $y$ is between $x,z$ if the cuts in $x$, $z$ divide different elements of the coarse partition of $\mathcal{E}$ induced by $y$, and $y$ is not equal to $x,z$.
\begin{pret5}
This relation defines a pretree.
\begin{proof}
Let
\begin{IEEEeqnarray*}{rCl}
\mathcal{E}&=&\pnm{x}{1}\sqcup\ldots\sqcup\pnm{x}{n_x}\\
\mathcal{E}&=&\pnm{y}{1}\sqcup\ldots\sqcup\pnm{y}{n_y}\\
\mathcal{E}&=&\pnm{z}{1}\sqcup\ldots\sqcup\pnm{z}{n_z}
\end{IEEEeqnarray*}
be the coarse partitions of the ends of $\Gamma$ induced by $x,y,z$.
First we check that the definition makes sense, i.e. given $x,y\in\mathcal{P}$ there are unique $i,j$ with 
\begin{equation*}
\coprod_{k\neq i} \pnm{x}{k}\subseteq \pnm{y}{j}
\end{equation*}
If one of $x,y$ is an equivalence class of totally isolated cuts, then lemmas \ref{ticc} and \ref{titi} yield this. Suppose both are \#-classes $Q,R$. By lemma \ref{cccc}, given $K\in R$ there is \pn{Q}{i}\ such that
\begin{equation*}
\pnm{K}{2}\subseteq\pnm{Q}{i}
\end{equation*} 
\pn{Q}{i}\ is contained in a \pn{\bar{Q}}{i}, so 
\begin{equation*}
\pnm{K}{2}\subseteq\pnm{\bar{Q}}{i}
\end{equation*} 
and furthermore this \pn{\bar{Q}}{i}\ is independent of the cut $K$ chosen. For each $j,j'$ we can find $K\in R$ with \pn{\bar{R}}{j}, \pn{\bar{R}}{j'} in different \pn{K}{k}\ since we are using the coarse partition; whence one of \pn{\bar{R}}{j}, \pn{\bar{R}}{j'} is contained in \pn{\bar{Q}}{i}. Hence all but one \pn{\bar{R}}{j}\ is contained in \pn{\bar{Q}}{i}\, i.e.
\begin{equation*}
\coprod_{k\neq i} \pnm{\bar{R}}{k}\subseteq \pnm{\bar{Q}}{i}
\end{equation*}
For part 1 of the definition of a pretree, note that if
\begin{IEEEeqnarray*}{rCl}
\coprod_{k\neq i_1} \pnm{x}{k}&\subseteq &\pnm{y}{j_1}\\
\coprod_{k\neq i_2} \pnm{x}{k}&\subseteq &\pnm{y}{j_2}
\end{IEEEeqnarray*}
with $j_1\neq j_2$ then since \pn{y}{j_1}, \pn{y}{j_2} are disjoint, $n_x=2=n_y$ and $x,y$ are equivalent cuts, hence equal as elements of $\mathcal{P}$. So $xyx$ does not hold.\\
2 is trivial. For 3, after relabelling we have
\begin{IEEEeqnarray*}{rCl}
\coprod_{k\neq 1} \pnm{y}{k}&\subseteq &\pnm{x}{1}\\
\coprod_{k\neq 1} \pnm{x}{k}&\subseteq &\pnm{y}{1}\\
\coprod_{k\neq 2} \pnm{y}{k}&\subseteq &\pnm{z}{1}\\
\coprod_{k\neq 1} \pnm{z}{k}&\subseteq &\pnm{y}{2}
\end{IEEEeqnarray*}
Then $\pnm{y}{1}\cap\pnm{y}{2}=\emptyset$ implies $\pnm{x}{i}\cap\pnm{z}{j}=\emptyset$ for $i,j\neq 1$. Hence 
\begin{IEEEeqnarray*}{rCl}
\coprod_{k\neq 1} \pnm{x}{k}&\subseteq &\pnm{z}{1}
\end{IEEEeqnarray*}
so $x,y$ divide the same \pn{z}{k}, hence $xzy$ does not hold.\\
For 4, suppose that $xzy$ so that 
\begin{IEEEeqnarray*}{rCl}
\coprod_{k\neq 1} \pnm{z}{k}&\subseteq &\pnm{x}{1}\\
\coprod_{k\neq 2} \pnm{z}{k}&\subseteq &\pnm{y}{1}
\end{IEEEeqnarray*}
i.e. $x$ divides \pn{z}{1}, $y$ divides \pn{z}{2}. If $w\neq z$ then $w$ divides a unique \pn{z}{i}. If $i=1$ then $yzw$. If not, then $xzw$.
\end{proof}
\end{pret5}
We recall the vertex version of Menger's Theorem (see for instance \cite{bondy2008graph}, thm 9.1, page 208):
\begin{menger}
Let $\Gamma$ be a graph and $a, b$ be vertices of $\Gamma$. Then the minimum size of a vertex cut separating $a, b$ is equal to the maximum number of vertex-independent simple paths joining $a,b$. 
\end{menger}
\begin{pret6}\label{discrete}
This pretree is discrete.
\begin{proof}
Let $K,L,M$ be mincuts with $M$ between $K,L$. Elements of $\mathcal{P}$ are of course not cuts; take a representative cut of any equivalence class or an appropriate corner cut of a \#-class. By lemma \ref{nestfin} only finitely many cuts are not nested with both $K,L$, so we need only consider the case when $M$ is nested with both. Let $C_i, D_i, E_i$ denote components of $\Gamma - K$, $\Gamma - L$, $\Gamma - M$ respectively. We have that $K$ is nested with $M$, so (after relabelling if necessary) $C_1\subseteq E_1$, and similarly $D_1\subseteq E_2$. $E_1\neq E_2$ as $M$ is between $K,L$. By the remarks following \ref{nestdef}, $E_1$ is contained in a component $D_i$, whence $K, L$ are nested.\\
If we now form a new graph by collapsing both of $K,L$ to a single vertex and apply Menger's Theorem in this graph, we obtain $n$ vertex-independent paths from $K$ to $L$, where $n$ is the cardinality of a mincut. In the case when $K,L$ are not disjoint then some of these paths collapse into points. The cut $M$ must intersect each of these paths as it separates $K, L$, and $|M|=n$, so $M$ is contained in the union of these paths. Then there are only finitely many choices for $M$.\\
If we took different choices for $K,L$ the only additional choices for $M$ would be equivalent in $\mathcal{P}$ to some already considered. So the pretree is discrete. 
\end{proof}
\end{pret6}
We now have a discrete pretree $\mathcal{P}$, which as discussed above gives us a tree encoding the mincuts of $\Gamma$ and how they interact with the ends of $\Gamma$.

%% file: succulents.tex
We have now obtained a tree encoding the cuts of the graph, with \#-classes collapsed down to points. We now seek to reintroduce the cyclic structure of these in order to obtain the final ``cactus" theorem. We will not be able to use cactus graphs as such; these work well for encoding edge cuts, but cannot represent a vertex cut yielding several components. We will therefore use a slightly more general structure which, for the sake of a horticultural joke, we call succulents.
\begin{succdef}\label{sucdef}
A \emph{succulent} is a connected graph built up from cycles (including possibly 2-cycles, consisting of two vertices joined by a double edge) in the following manner. Two cycles may be joined together either at a single vertex or along a single edge. The construction is tree-like in the sense that if we have a ``cycle of cycles" $C_1,\ldots,C_n$ with $C_i$ attached to $C_{i+1}$ (mod $n$) then all the $C_i$ share a common edge/vertex. The analogous property in a tree is that if we have a cycle of edges with each attached to the next one, they all meet at a common vertex. An \emph{end vertex} of a succulent is one contained in only one cycle; a vertex of a succulent is an end vertex if it has at most two edges adjacent to it.
\end{succdef}
\begin{figure}[htp]
\centering
\includegraphics[width=0.5\textwidth]{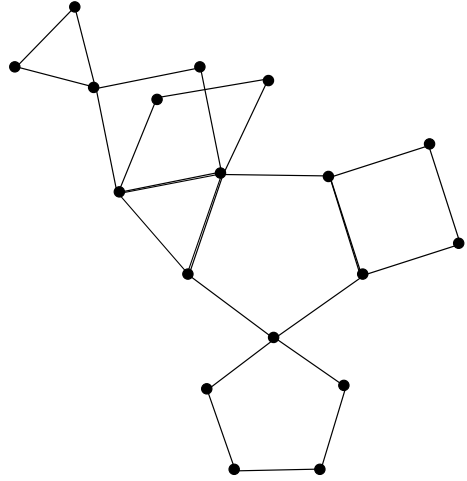}
\caption{A succulent.}
\end{figure}
We now construct a succulent encoding the mincuts of $\Gamma$. We already have the tree $T$ from the previous section whose vertices are (equivalence classes of) totally isolated cuts and \#-classes joined together via ``star'' vertices. There is at most one star for each corner cut of a \#-class. If there is no corresponding star, then the components split off by this cut are not further subdivided by mincuts.\\
Before moving on further we note that totally isolated cuts can be represented by a degenerate sort of ring, constructed by attaching traingles to a segment rather than to a cycle. So we can always talk about the anchors of a member of $\mathcal{P}$.\\
To form our succulent we replace each member of $\mathcal{P}$ by its associated ring. We must now consider how we connect these up; i.e. we need to consider the behaviour around each star vertex. Recall that if $Q$ is a \hashc attached to some star vertex, all members of $\mathcal{P}$ divide the same member of the (coarse) partition of the ends of $\Gamma$ corresponding to $Q$, \pn{\bar{Q}}{1} say, and that there is a corner cut of $Q$ separating \pn{\bar{Q}}{1} from the rest of the ends.\\
Suppose $Q,R$ are \#-classes adjacent to the same star vertex, so that there is no member of $\mathcal{P}$ between them. Let $K,L$ be the corresponding corner cuts and \pn{\bar{Q}}{1}, \pn{\bar{R}}{1} the members of the coarse partition. If both $K,L$ are B-cuts then each of \pn{\bar{Q}}{1}, \pn{\bar{R}}{1} comprises only one \pn{Q}{i}, \pn{R}{i} and there is only one member of the fine partitions divided by the other \#-class. We join these classes up by identifying the appropriate anchors. If there are no other elements of $\mathcal{P}$ joined to this star vertex then $K,L$ are equivalent so we could further simplify things by removing the anchors and joining the cycles for $Q,R$ together directly.\\
\begin{figure}[htp]
\centering
\begin{subfigure}[b]{0.3\textwidth}
\centering
\includegraphics[width=\textwidth]{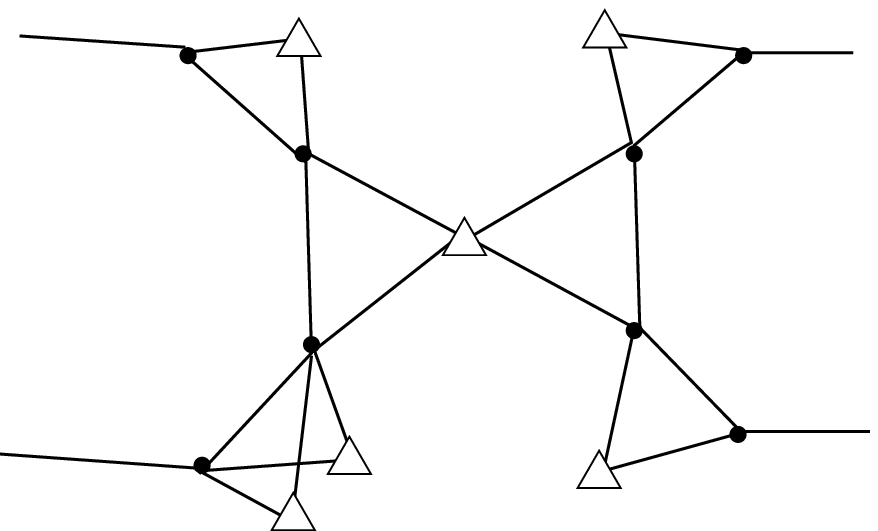}
\caption{If both corner cuts are B-cuts.}
\end{subfigure}%
~ 
\begin{subfigure}[b]{0.3\textwidth}
\centering
\includegraphics[width=\textwidth]{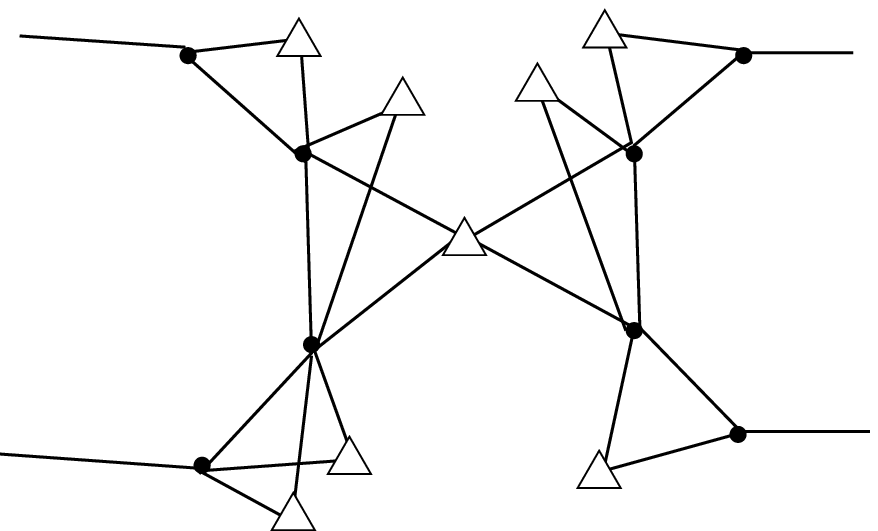}
\caption{If the corner cuts are nested and not equal.}
\end{subfigure}
~ 
\begin{subfigure}[b]{0.2\textwidth}
\centering
\includegraphics[width=\textwidth]{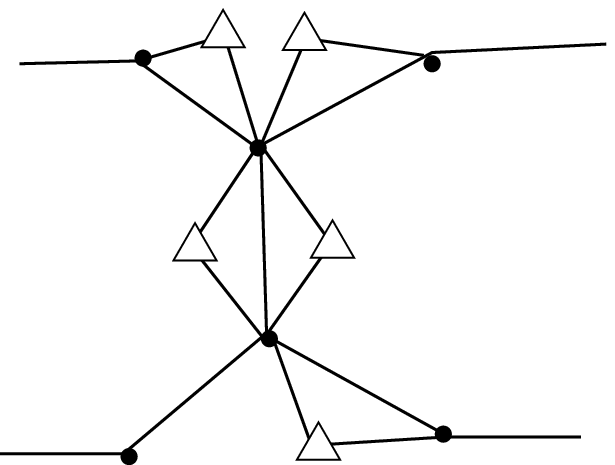}
\caption{If they are equal.}
\end{subfigure}
\caption{Diagrams illustrating how we connect up rings around a star vertex.}
\end{figure}
If one of $K,L$ is not a B-cut then the two cuts are nested. Then either they are equal or all components except one of $\Gamma - L$ are contained in the same component of $\Gamma - K$ and vice versa. In the latter case, there is only one  member of the fine partitions divided by the other \#-class, so again we can represent this by identifying the appropriate anchors. If the corner cuts are equal, then we glue together the rings via the corner cuts. Of the anchors attached to each, one represents the other \hashc and we simply delete this; the other anchors come in pairs, each representing the same set of ends but originating from different rings; we identify these together so we don't get redundancy. Then to produce our succulent we first glue together those rings sharing a corner cut, and then attach the other members of $\mathcal{P}$ adjacent to this star by identifying the appropriate anchors (for totally isolated cuts we simply note that the coarse and fine partitions coincide so there will be an obvious anchor to use and we have none of the issues above).\\
We must now show that this is a true succulent, that is that we still have a tree-like structure. We inherit much of the tree-like nature from $T$; we only need check that no ``cycles of cycles" form from the identifications made between rings all adjacent to the same star vertex. We will proceed by contradiction, supposing we have a shortest cycle of cycles $\mathcal{C}$ preventing our graph being a succulent. We can place limitations on which constituent cycles of the graph can be present in $\mathcal{C}$. First the cycles on which our rings are based do not appear. This is because any two cycles meeting one of these in the same \pn{\bar{Q}}{i} intersect along an edge. So $\mathcal{C}$ consists of the triangles which contain anchors; these can be joined together either at an anchor or along the opposite edge. Because our cycle is shortest, we alternate between joins along edges and at anchors. Hence our cycle has at least four members. Let $T_1,\ldots, T_4,\ldots$ be the triangles in $\mathcal{C}$ with $T_1, T_2$ meeting at an anchor, $T_2, T_3$ at an edge and so on. By construction the points at the bases of the $T_i$ represent cuts $K_i$ of $\Gamma$ partitioning the ends of $\Gamma$, and after suitable labelling we have 
\begin{eqnarray*}
\coprod_{i\neq 1} \pnm{K_1}{1}\subseteq \pnm{K_2}{1} = \pnm{K_3}{1}\\
\coprod_{i\neq 1} \pnm{K_4}{1}\subseteq \pnm{K_2}{2} = \pnm{K_3}{2}\\
\end{eqnarray*}
whence
\begin{IEEEeqnarray*}{rCl}
\coprod_{i\neq 1} \pnm{K_1}{1}&\subseteq&\pnm{K_4}{1}\\
&\subseteq& \pnm{K_6}{1}\\
&\ldots &
\end{IEEEeqnarray*}
But $\mathcal{C}$ is a cycle, so we eventually come back to the start whence all the inequalities become equalities, $K_1$ becomes a B-cut and
\begin{equation*}
\pnm{K_1}{2}=\pnm{K_2}{1}=\pnm{K_4}{1}=\ldots
\end{equation*}
we could have started at any other point, so the other $K_i$ are also B-cuts and all of them are equivalent. Then $\mathcal{C}$ becomes trivial and we have indeed constructed a succulent.\\
We have now proved most of:
\begin{succthm}\label{sucthm}
Let $\Gamma$ be a connected graph such that there are vertex end cuts of $\Gamma$ with finite cardinality. There is a succulent $\mathcal{S}$ with the following properties:
\begin{enumerate}
\item There is a subset $A$ of vertices of $\mathcal{S}$ called the \emph{anchors} of $\mathcal{S}$. If two anchors are adjacent, one of them is an end vertex of the graph. Every vertex of $\mathcal{S}$ not in $A$ is adjacent to an anchor. We define an anchor cut of $\mathcal{S}$ to be a vertex cut containing no anchors which separates some anchors of $\mathcal{S}$. We say anchor cuts are equivalent if they partition $A$ in the same way.
\item There is an onto map $f$ from the ends of $\Gamma$ to the union of the ends of $\mathcal{S}$ with the end vertices of $\mathcal{S}$ which are anchors.
\item There is a bijective map $g$ from equivalence classes of minimal end cuts of $\Gamma$ to equivalence classes of minimal anchor cuts of $\mathcal{S}$ such that ends $e_1, e_2$ of $\Gamma$ are separated by $[K]$ if and only if $f(e_1), f(e_2)$ are separated by $g([K])$.
\item Any automorphism of $\Gamma$ induces an automorphism of $\mathcal{S}$.
\end{enumerate}
\begin{proof}
We already have a succulent containing a representative of each mincut, i.e. we already have the map $g$. We now discuss how we modify the succulent to define the map $f$ of the ends of $\Gamma$. Some issues arise because there may be ends of $\Gamma$ which are distinguished from each other only by non-minimal cuts. If such ends exist, we will treat them as a single end for the present section, i.e. we will map them all to the same place using $f$.
Let $\epsilon$ be an end of $\Gamma$. If there is a mincut $K$ such that this end is the sole element of one of the sets \pn{K}{i}, then this mincut appears somewhere in the succulent either as a corner cut of a \hashc or as a totally isolated cut and there is an end anchor of the succulent corresponding to this \pn{K}{i}; define $f(\epsilon)$ to be this anchor.\\
If not, there may be a sequence of $x_i \in \mathcal{P}$, with 
\begin{equation*}
\pnm{x_1}{1}\supseteq\pnm{x_2}{1}\supseteq\ldots\ni\epsilon
\end{equation*}
This defines a ray in the tree $T$ associated to $\mathcal{P}$, hence an end of that tree. There is a unique such end, since $T$ is a tree so two ends can be separated using a single point, which we make take to be some $y\in\mathcal{P}$. But there is only one \pn{y}{i} containing $\epsilon$, so only one end will do. So we have an end of $T$, hence of the succulent, associated with $\epsilon$; this is where we will map $\epsilon$.\\
The remaining cases will correspond to ends which can only be split off by non-minimal cuts, which are not associated to some end of the tree $T$. To fit these into our succulent, we will essentially pretend that they can be split off by a mincut; we will add an element to $\mathcal{P}$ for each such end, inducing a partition
\begin{equation*}
\mathcal{E}=\{\epsilon\}\cup(\mathcal{E}-\{\epsilon\})
\end{equation*}
This member of $\mathcal{P}$ will not be between any two members of $\mathcal{P}$; and it does not disrupt the discreteness of $\mathcal{P}$ because an infinite set of betweeness in $\mathcal{P}$ would induce a descending sequence of partitions as above, so we would already have dealt with this end. So we have added an end vertex to the tree $T$. When constructing the succulent, the extra member of $\mathcal{P}$ will be modelled as two anchors joined by a double edge, one of which becomes attached to a relevant anchor in the succulent. The other anchor is an end anchor, which we define to be $f(\epsilon)$.\\
We now have the map $f$, which by construction interacts with $g$ in the way stated; note that the extra anchors added in the third step above are never split off by an anchor cut of $\mathcal{S}$.\\
To see that $f$ is onto, we note that any end anchors of $\mathcal{S}$ arise either as $f(\epsilon)$ in the third case above, or as part of a ring, where they correspond to some member of a partition of $\mathcal{E}$, whose members will be mapped there. Any ends of $\mathcal{S}$ arise from ends of $T$, hence from sequences of members of $\mathcal{P}$. From the vertices in the relevant cuts we can construct a ray in $\Gamma$, giving an end that will be mapped to the end of $\mathcal{S}$.\\
(4) arises since an automorphism of $\Gamma$ induces corresponding automorphisms of the cuts and ends of $\Gamma$, preserving crossings, nestings, equivalences; in short, all the information used to construct $\mathcal{S}$.
\end{proof}
\end{succthm}
We make some remarks about the theorem. In (3) we must say equivalence classes of cuts of $\mathcal{S}$ because we may have equivalent distinct cuts of $\mathcal{S}$; these arise if there are quasi-equivalent, non-equivalent half-cuts in a \hashc whence there will be some equivalent cuts contained in the relevant ring; but this is not really a concern.\\
If we wish to obtain a graph in which we do not have to exclude anchors from cuts, we can replace each anchor with a 3-vertex and treat these as ends, so that the anchor cuts in the theorem become bona fide mincuts of the resulting graph $\mathcal{S'}$.\\
If we collapse the extra end anchors we added in the proof above onto the adjacent anchors, then we obtain a variant theorem:
\begin{succthm2}\label{sucthm2}
Let $\Gamma$ be a connected graph such that there are vertex end cuts of $\Gamma$ with finite cardinality. There is a succulent $\mathcal{S}$ with the following properties:
\begin{enumerate}
\item There is a subset $A$ of vertices of $\mathcal{S}$ called the \emph{anchors} of $\mathcal{S}$. No two anchors are adjacent, and every vertex of $\mathcal{S}$ not in $A$ is adjacent to an anchor. We define an anchor cut of $\mathcal{S}$ to be a vertex cut containing no anchors which separates some anchors of $\mathcal{S}$. We say anchor cuts are equivalent if they partition $A$ in the same way.
\item There is a map $f$ from the ends of $\Gamma$ to the union of the ends of $\mathcal{S}$ with the anchors of $\mathcal{S}$.
\item There is a bijective map $g$ from equivalence classes of minimal end cuts of $\Gamma$ to equivalence classes of minimal anchor cuts of $\mathcal{S}$ such that ends $e_1, e_2$ of $\Gamma$ are separated by $[K]$ if and only if $f(e_1), f(e_2)$ are separated by $g([K])$.
\item Any automorphism of $\Gamma$ induces an automorphism of $\mathcal{S}$.
\end{enumerate}
\end{succthm2}
Consider a finite graph $\Gamma$. We call a set $J$ of vertices of a graph $\Gamma$ \emph{$n$-inseparable} if $|J|\geq n+1$ and for any set $K$ of vertices with $|K|\leq n$, $J$ is contained in a single component of $\Gamma - K$. Let $\kappa$ be the smallest integer for which there $\kappa$-inseparable sets $J_1$, $J_2$ and a vertex cut $K$ with $|K|=\kappa$ and $J_1, J_2$ in different components of $\Gamma - K$. We can consider the maximal $\kappa$-inseperable sets of $\Gamma$ as ends of the graph; or attach a sequence of $(\kappa + 1)$-vertices to each to turn them into a bona fide end. The inseparability conditions ensure that this does not affect the cuts of $\Gamma$ of size $\kappa$. Then the size $\kappa$ vertex cuts separating two inseperable sets become minimal end cuts of our graph, so we can obtain a succulent theorem for them:
\begin{succthm3}
Let $\Gamma$ be a finite connected graph such that there exists $\kappa$  for which there $\kappa$-inseparable sets $J_1$, $J_2$ and a vertex cut $K$ with $|K|=\kappa$ and $J_1, J_2$ in different components of $\Gamma - K$, and take the minimal such $\kappa$. There is a succulent $\mathcal{S}$ with the following properties:
\begin{enumerate}
\item There is a subset $A$ of vertices of $\mathcal{S}$ called the \emph{anchors} of $\mathcal{S}$. No two anchors are adjacent, and every vertex of $\mathcal{S}$ not in $A$ is adjacent to an anchor. We define an anchor cut of $\mathcal{S}$ to be a vertex cut containing no anchors which separates some anchors of $\mathcal{S}$. We say anchor cuts are equivalent if they partition $A$ in the same way.
\item There is a map $f$ from the $\kappa$-inseparable sets of $\Gamma$ to the anchors of $\mathcal{S}$.
\item There is a bijective map $g$ from equivalence classes of minimal cuts of $\Gamma$ separating $\kappa$-inseparable sets to equivalence classes of minimal anchor cuts of $\mathcal{S}$ such that $\kappa$-inseparable sets $J_1, J_2$ of $\Gamma$ are separated by $[K]$ if and only if $f(J_1), f(J_2)$ are separated by $g([K])$.
\item Any automorphism of $\Gamma$ induces an automorphism of $\mathcal{S}$.
\end{enumerate}
\end{succthm3}
Tutte \cite{tutte1984graph} produced structure trees for the cases $\kappa = 1,2$, which Dunwoody and Kr\"on \cite{dunwoody2009vertex} then extended to higher $\kappa$. These trees were based on `optimally nested' cuts in the language of \cite{dunwoody2009vertex}, which in this case means A-cuts. Roughly speaking, the trees consist of the totally isolated cuts and corner cuts of our succulents, together with `blocks' which are not decomposed by the cuts in question; these include the maximal inseparable sets, and also sets broken up by cuts which are not optimally nested; these sets correspond to the \#-classes. The structure trees can then be obtained from our succulents by replacing each ring with a star with one central vertex and one vertex joined to it for each corner cut. So these earlier results also follow from our work.

%% file: applications.tex
First we note that our work yields a proof of Stallings's Theorem, based on the Bass-Serre theory of groups acting on trees (see \cite{serre1980trees}).
\begin{stallings}
Let $G$ be a finitely generated group acting transitively on a graph $\Gamma$ with more than two ends. Then $G$ can be expressed as an amalgam $G=A*_F B$ or an HNN extension $G=A*_F$ where $F$ has a finite index subgroup which is the stabilizer of a vertex of $\Gamma$.
\begin{proof}
From the pretree $\mathcal{P}$ we obtain a tree $T$ on which $G$ acts. The tree $T$ is non-trivial; the action is transitive and $\Gamma$ has more than two ends so there are infinitely many ends and many inequivalent cuts. The action is without inversion since $T$ is bipartite, formed of star vertices and elements of $\mathcal{P}$. Then $G$ is isomorphic to the fundamental group of a certain graph of groups; $G$ is finitely generated so this graph is finite. The action is non-trivial as $G$  acts transitively on $\Gamma$, so it follows that $G$ splits over the stabilizer of an edge of $T$. An element fixing an edge of $T$ fixes the adjacent element of $\mathcal{P}$, hence fixes either a \hashc or an equivalence class of totally isolated cuts. A \hashc contains finitely many vertices; and the transitivity of the action implies that there can only be finitely many cuts in each equivalence class, since we can find two cuts between which every cut of the class lies, and then apply the methods of \ref{discrete}. The result follows.
\end{proof}
\end{stallings}
Stallings's original theorem covers the two-ended case as well, but our tree is trivial here. The two-ended case can be covered by more elementary means however.\\
We now discuss how earlier cactus theorems concerning edge cuts follow from ours. We turn a question about edge end cuts into a question about vertex end cuts as follows. First replace the graph $\Gamma$ with its barycentric subdivision $\Gamma^b$. This is defined as follows: 
\begin{equation*}
V(\Gamma^b)=V(\Gamma)\cup E(\Gamma)
\end{equation*}
\begin{equation*}
E(\Gamma^b)=\{(v,e): v\in V(\Gamma), e\in E(\Gamma), v\text{ an endpoint of }e\}
\end{equation*}
If the cardinality of a minimal edge end cut of $\Gamma$ is $n$, then we now `thicken up' each vertex of $\Gamma^b$ that was a vertex of $\Gamma$ by replacing it with an $(n+1)$-vertex (see definition \ref{rindef}) to obtain a graph $\Gamma^*$. In this way an edge cut of $\Gamma$ separating some ends of $\Gamma$ corresponds precisely with a vertex cut of $\Gamma^*$ of the same cardinality. In $\Gamma^*$, because all the vertex cuts are essentially edge cuts, all of the minimal vertex cuts of $\Gamma^*$ split the graph into precisely two pieces, each containing an end. So we do not need to remove slices from the graph, and all cuts are B-cuts. It follows that quasi-equivalent half-cuts are equivalent, and each ring becomes simple enough to be replaced by a cycle, in which the anchors become the vertices and the other vertices become the edges. Our succulent from theorem \ref{sucthm2} can then be replaced with a cactus, so we have the cactus theorem for edge end cuts \cite{papaedgecuts}.
\begin{cactends}
Let $\Gamma$ be a connected graph such that there are edge end cuts of $\Gamma$ with finite cardinality. There is a cactus $\mathcal{C}$ with the following properties:
\begin{enumerate}
\item There is a map $f$ from the ends of $\Gamma$ to the union of the ends of $\mathcal{C}$ with the vertices of $\mathcal{C}$.
\item There is a bijective map $g$ from equivalence classes of minimal end cuts of $\Gamma$ to minimal edge cuts of $\mathcal{C}$ such that ends $e_1, e_2$ of $\Gamma$ are separated by $[K]$ if and only if $f(e_1), f(e_2)$ are separated by $g([K])$.
\item Any automorphism of $\Gamma$ induces an automorphism of $\mathcal{C}$.
\end{enumerate}
\end{cactends}
To deal with the classical cactus theorem for edge cuts of finite graphs, we proceed as before to get the graph $\Gamma^*$. Then to each $(n+1)$-vertex we attach an infinite chain of $(n+1)$-vertices, so that a vertex in the original graph $\Gamma$ becomes a de facto end of our new graph. `Equivalent cuts' in this graph correspond to the same cut of the original graph. Once again the succulent can be replaced with a cactus, so we have the cactus theorem of Dinits-Karzanov-Lomonosov \cite{dinits1991structure}.
\begin{cactedge}
Let $\Gamma$ be a connected finite graph. There is a cactus $\mathcal{C}$ with the following properties:
\begin{enumerate}
\item There is a map $f$ from the vertices of $\Gamma$ to the vertices of $\mathcal{C}$.
\item There is a bijective map $g$ from equivalence classes of minimal edge cuts of $\Gamma$ to minimal edge cuts of $\mathcal{C}$ such that vertices $v_1, v_2$ of $\Gamma$ are separated by $[K]$ if and only if $f(v_1), f(v_2)$ are separated by $g([K])$.
\item Any automorphism of $\Gamma$ induces an automorphism of $\mathcal{C}$.
\end{enumerate}
\end{cactedge}